\documentclass[12pt,a4paper]{article}

\usepackage[left=2cm,right=2cm,    top=2cm,bottom=2cm,bindingoffset=0cm]{geometry}

\usepackage{subfig}
\usepackage{flushend}
\usepackage{indentfirst}
\usepackage{graphics}
\usepackage{amsmath}
\usepackage{epstopdf}
\usepackage{float}
\usepackage{amssymb}
\usepackage[font=footnotesize,labelfont=bf]{caption}
\usepackage[labelsep=period]{caption}

\newcommand{\be}{\begin{equation}}
\newcommand{\ee}{\end{equation}}

\newtheorem{thm}{Theorem}[section]
\newtheorem{lem}{Lemma}[section]
\newtheorem{nas}{Corollary}[section]

\newtheorem{zau}{Remark}[section]

\newtheorem{ozn}{Definition}[section]

\begin{document}

\title{\textbf{On Minimax Estimation Problems for Periodically Correlated Stochastic Processes }}

\date{}

\maketitle

\noindent Columbia International Publishing\\
Contemporary Mathematics and Statistics
(2014) Vol. 2 No. 1 pp. 1-24\\
doi:10.7726/cms.2014.1001

\vspace{20pt}

\author{\textbf{Iryna Dubovets'ka}, \textbf{Mykhailo Moklyachuk}$^{*}$, \\\\
 {Department of Probability Theory, Statistics and Actuarial
Mathematics, \\
Taras Shevchenko National University of Kyiv, Kyiv 01601, Ukraine}\\
$^{*}$Corresponding Author: Moklyachuk@gmail.com}\\\\\\

\noindent \textbf{Abstract.} \hspace{2pt}
The aim of this article is to overview the problem of mean square optimal estimation of linear functionals which depend on unknown values of periodically correlated stochastic process. Estimates are based on observations of this process and noise. These problems are investigated under conditions of spectral certainty and spectral uncertainty. Formulas for calculating the main characteristics (spectral characteristic, mean square error) of the optimal linear estimates of the functionals are proposed. The least favorable spectral densities and the minimax-robust spectral characteristics of optimal estimates of the functionals are presented for given sets of admissible spectral densities.\\

\noindent \textbf{Keywords:} \hspace{2pt}
 Periodically Correlated Process; Spectral Characteristic; Mean Square Error; Minimax Estimate; Minimax-Robust Spectral Characteristic\\

\noindent \textbf{Mathematics Subject Classification:} \hspace{2pt}
 Primary: 60G10, 60G25, 60G35, Secondary: 62M20, 93E10, 93E11

\section{{Introduction}}

Periodically correlated processes are those signals whose statistics vary almost periodically, and they are present in numerous physical and man-made processes. A comprehensive listing most of the existing references up to the year 2005 on periodically correlated processes and their applications was proposed by Serpedin et al. (2005). See also a review by Antoni (2009). For more details see survey paper by Gardner (1994) and book by Hurd and Miamee (2007). Note, that most of authors investigate properties of periodically correlated sequences while only few publications deal with investigation of periodically correlated processes. Note also, that in the literature periodically correlated processes are named in multiple different ways such as cyclostationary, periodically nonstationary or cyclic correlated processes.

Periodically correlated processes can be defined as stochastic processes with a periodic structure. In papers by Gladyshev (1961, 1963) investigation of periodically correlated processes was started. Analysis of properties of correlation function and representations of periodically correlated processes were presented. Relations between periodically correlated processes and stationary processes were investigated by Makagon (1999a, 2001). Relations of periodically correlated sequences with simpler stochastic sequences are proposed by Makagon (1999b, 2011), Makagon and Miamee (2013), Hurd and Miamee (2007).

Methods of solution of problems of estimation of unknown values of stationary stochastic processes (extrapolation, interpolation and filtering problems) were developed by Wiener (1966), Yaglom (1987), Kolmogorov (1992). Estimation problems for stationary vector sequences were investigated by Rozanov (1967). The proposed methods are based on the assumption that spectral densities of processes are exactly known. In practice, however, it is impossible to have complete information on the spectral density in most cases. To solve the problem one finds parametric or nonparametric estimates of the unknown spectral density or selects a density by other reasoning. Then the classical estimation method is applied provided that the estimated or selected density is the true one. This procedure can result in significant increasing of the value of error as Vastola and Poor (1983) have demonstrated with the help of some examples. This is a reason to search estimates which are optimal for all densities from a certain class of admissible spectral densities. These estimates are called minimax since they minimize the maximal value of the error. A survey of results in minimax (robust) methods of data processing can be found in the paper by Kassam and Poor (1985). The paper by Ulf Grenander (1957) should be marked as the first one where the minimax extrapolation problem for stationary processes was formulated and solved. Franke and Poor (1984), Franke (1984, 1985) investigated the minimax extrapolation and filtering problems for stationary sequences with the help of convex optimization methods. This approach makes it possible to find equations that determine the least favorable spectral densities for various classes of admissible densities. For more details see, for example, books by Moklyachuk (2008), Moklyachuk and Masyutka (2012). In papers by Moklyachuk (1994-2008) the minimax approach was applied to extrapolation, interpolation and filtering problems for functionals which depend on the unknown values of stationary processes and sequences. Methods of solution the minimax-robust estimation problems for vector-valued stationary sequences and processes were developed by Moklyachuk and Masyutka (2006-2011). Luz and Moklyachuk (2012-2013) investigated the minimax estimation problems for linear functionals which depends on unknown values of stochastic sequence with stationary th increments. Minimax estimation problems for linear functionals which depend on the unknown values of periodically correlated sequences and processes were studied in works by Dubovetska, Masyutka, Moklyachuk (2011 -2013). These problems are natural generalization of the extrapolation, interpolation and filtering problems for functionals which depend on the unknown values of stationary processes and sequences.

In this article we consider the problem of optimal linear estimation of functionals
$$A_{N}^{i}\zeta =\int_{0}^{(N+1)T}{a}(t)\zeta (t)dt,\quad {{A}^{e}}\zeta =\int_{0}^{\infty }{a(t)\zeta (t)\,dt},\quad A_{{}}^{f}\zeta =\int_{0}^{\infty }{a}(t)\zeta (-t)dt$$
which depend on the unknown values of a periodically correlated process $\zeta (t)$ based on observations of the process $\zeta (t)+\theta (t)$ at points of time  $t\in \mathbb{R}\backslash [0,(N+1)T]$ (for estimation  $A_{N}^{i}\zeta $), at points of time $t<0$ (for estimation ${{A}^{e}}\zeta $), at points of time $t\le 0$ (for estimation ${{A}^{f}}\zeta $). Here $\theta (t)$ is an uncorrelated with $\zeta (t)$ periodically correlated stochastic process. We propose a transition procedure from continuous periodically correlated stochastic processes $\{\zeta (t),\,t\in \mathbb{R}\}$ and $\{\theta (t),\,t\in \mathbb{R}\}$ to the corresponding infinite dimensional vector-valued stationary sequences $\{{{\zeta }_{j}},j\in \mathbb{Z}\}$ and $\{{{\theta }_{j}},j\in \mathbb{Z}\}$  which allows us to reduce the estimation problems for continuous periodically correlated stochastic processes to the corresponding problems for stationary vector-valued sequences.

Formulas for calculating the mean square errors and spectral characteristics of the optimal linear estimates of the corresponding functionals are proposed in the case of spectral certainty where spectral densities of generated stationary sequences $\{{{\zeta }_{j}},j\in \mathbb{Z}\}$ and $\{{{\theta }_{j}},j\in \mathbb{Z}\}$ are exactly known. The least favorable spectral densities and the minimax-robust spectral characteristics of the optimal linear estimates are found in the case of spectral uncertainty where spectral densities are not exactly known, but concrete classes of admissible densities are given. It is shown, for example, that one-sided moving average sequence gives the greatest value of the mean square error of the optimal estimate of the functional${{A}^{e}}\zeta $.

\section{{Periodically correlated continuous processes and generated vector-valued stationary sequences}}

\begin{ozn} \label{def2.1}
 (Gladyshev, 1963) Mean square continuous stochastic process $\zeta :\mathbb{R}\to H={{L}_{2}}(\Omega ,F,P)$, $E\zeta (t)=0$, is called periodically correlated (PC) with period $T$, if its correlation function $K(t,s)=E\zeta (t)\overline{\zeta (s)}$ for all $t,s\in \mathbb{R}$ and some fixed $T>0$ is such that
$K(t,s)=K(t+T,s+T).$
\end{ozn}
Let $\{\zeta (t),\,t\in \mathbb{R}\}$ and $\{\theta (t),\,t\in \mathbb{R}\}$ be uncorrelated PC stochastic processes with period $T$.
We construct the following sequences of stochastic functions
\begin{equation} \label{eq01}
\left\{ {{\zeta }_{j}}(u)=\zeta (u+jT),\,u\in [0,T),\,j\in \mathbb{Z} \right\},
\end{equation}
\begin{equation} \label{eq02}
\left\{ {{\theta }_{j}}(u)=\theta (u+jT),\,u\in [0,T),\,j\in \mathbb{Z} \right\}.
\end{equation}

Sequences (1) and (2) form ${{L}_{2}}([0,T);H)$-valued stationary sequences $\{{{\zeta }_{j}},\,j\in \mathbb{Z}\}$ and $\{{{\theta }_{j}},\,j\in \mathbb{Z}\}$,
respectively, with the correlation functions
     $${{B}_{\zeta }}(l,j)={{\left\langle {{\zeta }_{l}},{{\zeta }_{j}} \right\rangle }_{H}}=\int_{o}^{T}{E\,\zeta (u+lT)}\overline{\zeta (u+jT)}\,du=\int_{o}^{T}{{{K}_{\zeta }}(u+(l-j)T,u)}\,du={{B}_{\zeta }}(l-j),$$
          $${{B}_{\theta }}(l,j)={{\left\langle {{\theta }_{l}},{{\theta }_{j}} \right\rangle }_{H}}=\int_{o}^{T}{E\,\theta (u+lT)}\overline{\theta (u+jT)}\,du=\int_{o}^{T}{{{K}_{\theta }}(u+(l-j)T,u)}\,du={{B}_{\theta }}(l-j),$$
where ${{K}_{\zeta }}(t,s)=E\zeta (t)\overline{\zeta (s)},$ ${{K}_{\theta }}(t,s)=E\theta (t)\overline{\theta (s)}$ are correlation functions of PC processes $\zeta (t)$and $\theta (t)$. Consider in the space ${{L}_{2}}([0,T);\mathbb{R})$the following orthonormal basis
$$\{{{\tilde{e}}_{k}}=\frac{1}{\sqrt{T}}{{e}^{2\pi i\{{{(1)}^{-k}}\left[ {}^{k}/{}_{2} \right]\}u/T}},\,k=1,2,...\},\left\langle {{{\tilde{e}}}_{k}},{{{\tilde{e}}}_{j}} \right\rangle ={{\delta }_{jk}},
$$
Making use properties of this basis, stationary sequences $\{{{\zeta }_{j}},\,j\in \mathbb{Z}\}$ and $\{{{\theta }_{j}},\,j\in \mathbb{Z}\}$ can be represented in the forms
\begin{equation} \label{eq03}
{{\zeta }_{j}}=\sum\nolimits_{k=1}^{\infty }{{{\zeta }_{kj}}{{{\tilde{e}}}_{k}}},\quad {{\zeta }_{kj}}=\left\langle {{\zeta }_{j}},{{{\tilde{e}}}_{k}} \right\rangle =\frac{1}{\sqrt{T}}\int_{0}^{T}{{{\zeta }_{j}}(v){{e}^{-2\pi i\{{{(1)}^{-k}}\left[ {}^{k}/{}_{2} \right]\}v/T}}}\,dv,
 \end{equation}
\begin{equation} \label{eq04}
     {{\theta }_{j}}=\sum\nolimits_{k=1}^{\infty }{{{\theta }_{kj}}{{{\tilde{e}}}_{k}}},\quad {{\theta }_{kj}}=\left\langle {{\theta }_{j}},{{{\tilde{e}}}_{k}} \right\rangle =\frac{1}{\sqrt{T}}\int_{0}^{T}{{{\theta }_{j}}(v){{e}^{-2\pi i\{{{(1)}^{-k}}\left[ {}^{k}/{}_{2} \right]\}v/T}}}\,dv.
     \end{equation}
We call these sequences $\{{{\zeta }_{j}},\,j\in \mathbb{Z}\}$, $\{{{\theta }_{j}},\,j\in \mathbb{Z}\}$ and corresponding to them vector-valued sequences $\{{{\vec{\zeta }}_{j}}={{({{\zeta }_{kj}},\,k=1,2,...)}^{\top }},j\in \mathbb{Z}\}$, $\{{{\vec{\theta }}_{j}}={{({{\theta }_{kj}},\,k=1,2,...)}^{\top }},j\in \mathbb{Z}\},$
generated (by $\{\zeta (t),\,t\in \mathbb{R}\}$, $\{\theta (t),\,t\in \mathbb{R}\}$, respectively) vector-valued stationary sequences.

Components $\{{{\zeta }_{kj}},\,k=1,2,...\}$ and $\{{{\theta }_{kj}},\,k=1,2,...\}$ of generated stationary sequences $\{{{\zeta }_{j}},\,j\in \mathbb{Z}\}$ and $\{{{\theta }_{j}},\,j\in \mathbb{Z}\}$ are such that (Kallianpur and Mandrekar, 1971; Moklyachuk, 1981)
$$E{{\zeta }_{kj}}=0,\quad ||{{\zeta }_{j}}||_{H}^{2}=\sum\nolimits_{k=1}^{\infty }{E|{{\zeta }_{kj}}{{|}^{2}}}={{P}_{\zeta }}<\infty,\quad E{{\zeta }_{kl}}\overline{{{\zeta }_{nj}}}=\left\langle {{R}_{\zeta }}(l-j){{e}_{k}},{{e}_{n}} \right\rangle,
$$
 $$
 E{{\theta }_{kj}}=0,\quad ||{{\theta }_{j}}||_{H}^{2}=\sum\nolimits_{k=1}^{\infty }{E|{{\theta }_{kj}}{{|}^{2}}}={{P}_{\theta }}<\infty,\quad E{{\theta }_{kl}}\overline{{{\theta }_{nj}}}=\left\langle {{R}_{\theta }}(l-j){{e}_{k}},{{e}_{n}} \right\rangle,
 $$
where $\{{{e}_{k}},\,\,k=1,2,...\}$ is a basis of the space ${{\ell }_{2}}$.
Correlation functions ${{R}_{\zeta }}(j)$ and ${{R}_{\theta }}(j)$ of stationary sequences $\{{{\zeta }_{j}},\,j\in \mathbb{Z}\}$ and $\{{{\theta }_{j}},\,j\in \mathbb{Z}\}$ are correlation operator functions in the space ${{\ell }_{2}}$.
Correlation operators ${{R}_{\zeta }}(0)={{R}_{\zeta }}$, ${{R}_{\theta }}(0)={{R}_{\theta }}$ are kernel operators:
$$\sum\nolimits_{k=1}^{\infty }{\left\langle {{R}_{\zeta }}{{e}_{k}},{{e}_{k}} \right\rangle }=\,||{{\zeta }_{j}}||_{H}^{2}={{P}_{\zeta }},\quad \sum\nolimits_{k=1}^{\infty }{\left\langle {{R}_{\theta }}{{e}_{k}},{{e}_{k}} \right\rangle }=\,||{{\theta }_{j}}||_{H}^{2}={{P}_{\theta }}.$$
The stationary sequences $\{{{\zeta }_{j}},\,j\in \mathbb{Z}\}$ and $\{{{\theta }_{j}},\,j\in \mathbb{Z}\}$ have spectral densities
$$f(\lambda )=\{{{f}_{kn}}(\lambda )\}_{k,n=1}^{\infty },\quad g(\lambda )=\{{{g}_{kn}}(\lambda )\}_{k,n=1}^{\infty },$$
which are positive operator-valued functions in ${{\ell }_{2}}$ of the variable $\lambda \in [-\pi ,\pi )$, if theirs correlation functions
${{R}_{\zeta }}(j)$ and ${{R}_{\theta }}(j)$ can be represented in the form
$$\left\langle {{R}_{\zeta }}(j){{e}_{k}},{{e}_{n}} \right\rangle =\frac{1}{2\pi }\int_{-\pi }^{\pi }{{{e}^{ij\lambda }}\left\langle f(\lambda ){{e}_{k}},{{e}_{n}} \right\rangle }\,d\lambda ,$$
$$\left\langle {{R}_{\theta }}(j){{e}_{k}},{{e}_{n}} \right\rangle =\frac{1}{2\pi }\int_{-\pi }^{\pi }{{{e}^{ij\lambda }}\left\langle g(\lambda ){{e}_{k}},{{e}_{n}} \right\rangle }\,d\lambda ,$$ $k,n=1,2,....$
For almost all $\lambda \in [-\pi ,\pi )$ spectral densities $f(\lambda ),$ $g(\lambda )$ are kernel operators with the integrable kernel norms
$$\sum\nolimits_{k=1}^{\infty }{\frac{1}{2\pi }\int_{-\pi }^{\pi }{\left\langle f(\lambda ){{e}_{k}},{{e}_{k}} \right\rangle }\,d\lambda =\sum\nolimits_{k=1}^{\infty }{\,\left\langle {{R}_{\zeta }}{{e}_{k}},{{e}_{k}} \right\rangle }}=||{{\zeta }_{j}}||_{H}^{2}={{P}_{\zeta }},$$
$$\sum\nolimits_{k=1}^{\infty }{\frac{1}{2\pi }\int_{-\pi }^{\pi }{\left\langle g(\lambda ){{e}_{k}},{{e}_{k}} \right\rangle }\,d\lambda =\sum\nolimits_{k=1}^{\infty }{\,\left\langle {{R}_{\theta }}{{e}_{k}},{{e}_{k}} \right\rangle }}=||{{\theta }_{j}}||_{H}^{2}={{P}_{\theta }}.$$
We will use representations (3), (4) for finding solutions to the mean square estimation problems for continuous periodically correlated stochastic processes.

\section{{Hilbert space projection method of estimation of PC processes}}

\subsection{ Interpolation problem}

Consider the problem of optimal linear estimation of the functional
$$A_{N}^{i}\zeta =\int_{0}^{(N+1)T}{a}(t)\zeta (t)dt$$
which depends on the unknown values of the mean square continuous PC stochastic process $\zeta (t)$ based on observations of the process $\zeta (t)+\theta (t)$ at points of time $t\in \mathbb{R}\backslash [0,(N+1)T]$.
The noise process $\theta (t)$ is an uncorrelated with $\zeta (t)$ PC stochastic process.
To be sure that the functional $A_{N}^{i}\zeta $ is well defined we will suppose that the function $a(t),\,\,t\in {{\mathbb{R}}_{+}},$ satisfies the natural necessary condition $$\int_{0}^{(N+1)T}{|a(t)|\,dt}<\infty.$$
With the help of transformation (1) of the process $\zeta (t)$ we can represent the functional $A_{N}^{i}\zeta $ in the form
$$A_{N}^{i}\zeta =\int_{0}^{N}{a(t)\zeta (t)\,dt=}\sum\nolimits_{j=0}^{N}{\int_{0}^{T}{{{a}_{j}}(u){{\zeta }_{j}}(u)}\,du},$$
where ${{a}_{j}}(u)=a(u+jT),$ ${{\zeta }_{j}}(u)=\zeta (u+jT),$ $u\in [0,T).$
Making use the decomposition (3) of the generated stationary sequence $\{{{\zeta }_{j}},\,j\in \mathbb{Z}\}$ and solutions of the equation
$${{(-1)}^{k}}\left[ {k}/{2}\; \right]+{{(-1)}^{n}}\left[ {n}/{2}\; \right]=0$$
of two variables $(k,n)$, which are given by pairs $(1,1),$$(2l+1,2l)$ and $(2l,2l+1)$ for $l=2,3,...,$ the functional can be represented in the form
$$A_{N}^{i}\zeta =\sum\nolimits_{j=0}^{N}{\int_{0}^{T}{{{a}_{j}}(u){{\zeta }_{j}}(u)}\,du}=$$
$$=\sum\nolimits_{j=0}^{N}{\frac{1}{T}\int_{0}^{T}{\left( \sum\nolimits_{k=1}^{\infty }{{{a}_{kj}}\exp \left\{ {2\pi i\{{{(-1)}^{k}}\left[ {k}/{2}\; \right]\}u}/{T}\; \right\}} \right)}}
\times$$
$$\times\left( \sum\nolimits_{n=1}^{\infty }{{{a}_{nj}}\exp \left\{ {2\pi i\{{{(-1)}^{n}}\left[ {n}/{2}\; \right]\}u}/{T}\; \right\}} \right)\,du=$$
$$=\sum\nolimits_{j=0}^{N}{\sum\nolimits_{k=1}^{\infty }{\sum\nolimits_{n=1}^{\infty }{{{a}_{kj}}{{\zeta }_{kj}}\frac{1}{T}}}}\int_{0}^{T}{\exp \left\{ {2\pi i\{{{(-1)}^{k}}\left[ {k}/{2}\; \right]+{{(-1)}^{n}}\left[ {n}/{2}\; \right]\}u}/{T}\; \right\}}\,du=$$
$$=\sum\nolimits_{j=0}^{N}{\sum\nolimits_{k=1}^{\infty }{{{a}_{kj}}{{\zeta }_{kj}}=}}\sum\nolimits_{j=0}^{N}{\vec{a}_{j}^{\top }{{{\vec{\zeta }}}_{j}}},$$
where vectors $\,\,{{\vec{a}}_{j}}={{({{a}_{kj}},\,k=1,2,...)}^{\top }}={{({{a}_{1j}},{{a}_{3j,}}{{a}_{2j}},...,{{a}_{2k+1,j}},{{a}_{2k,j}},...)}^{\top }},$ ${{a}_{kj}}=\left\langle {{a}_{j}},{{{\tilde{e}}}_{k}} \right\rangle .$

Assume that vectors $\{{{\vec{a}}_{j}},\,j=0,1,...,N\}$ satisfy the following conditions
\begin{equation} \label{eq05}
||{{\vec{a}}_{j}}||<\infty ,\,\,\,\,\,||{{\vec{a}}_{j}}|{{|}^{2}}=\sum\nolimits_{k=1}^{\infty }{|{{a}_{kj}}{{|}^{2}}},\,\,\,\,j=0,1,...N.
\end{equation}
It follows from condition (5) that the functional $A_{N}^{i}\zeta $ has finite second moment.

Let spectral densities $f(\lambda )$ and $g(\lambda )$ of the generated stationary sequences $\{{{\zeta }_{j}},\,j\in \mathbb{Z}\}$ and $\{{{\theta }_{j}},\,j\in \mathbb{Z}\}$ be such that the minimality condition is satisfied
\begin{equation} \label{eq06}
\int_{-\pi }^{\pi }{Tr[{{(f(\lambda )+g(\lambda ))}^{-1}}]\,d\lambda <\infty \,.}
\end{equation}
The minimality condition (6) is necessary and sufficient in order that the error-free estimation of unknown values of the sequence $\{{{\zeta }_{j}}+{{\theta }_{j}},\,j\in \mathbb{Z}\}$ is impossible (Rozanov, 1967).

Denote by ${{L}_{2}}(f)$ the Hilbert space of vector-functions $b(\lambda )=\{{{b}_{k}}(\lambda )\}_{k=1}^{\infty }$, which are integrable with respect to the measure with density $f(\lambda )$
$$\int_{-\pi }^{\pi }{{{b}^{\top }}(\lambda )f(\lambda )}\overline{b(\lambda )}\,d\lambda =\int_{-\pi }^{\pi }{\sum\nolimits_{k,n=1}^{\infty }{{{b}_{k}}(\lambda ){{f}_{kn}}(\lambda )\overline{{{b}_{n}}(\lambda )}}}\,d\lambda <\infty.$$

Denote by $L_{2}^{\,N-}(f+g)$  the subspace of $L_{2}^{{}}(f+g)$ generated by vector-functions ${{e}^{ij\lambda }}{{\delta }_{k}},$ $j\in \mathbb{Z}\backslash \{0,1,...,N\},$ $k=1,2,...,$ where ${{\delta }_{kn}}$ is the Kronecker symbol: ${{\delta }_{kk}}=1$ and ${{\delta }_{kn}}=0$ for $k\ne n$.
Every estimate $\hat{A}_{N}^{i}\zeta $ of the functional $A_{N}^{i}\zeta $ based on observations of the process $\zeta (t)+\theta (t)$ at points $t\in \mathbb{R}\backslash [0,(N+1)T]$ is characterized by its spectral characteristic $h({{e}^{i\lambda }})\in L_{2}^{N\,-}(f+g)$ and the orthogonal stochastic measure ${{Z}^{\zeta +\theta }}(\Delta )=\{Z_{k}^{\zeta +\theta }(\Delta )\}_{k=1}^{\infty }$ of the sequence $\{{{\zeta }_{j}}+{{\theta }_{j}},j\in \mathbb{Z}\}$ and has the following form
\begin{equation} \label{eq07}
\hat{A}_{N}^{i}\zeta =\int_{-\pi }^{\pi }{{{h}^{\top }}({{e}^{i\lambda }})}({{Z}^{\zeta +\theta }}(\,d\lambda ))=\int_{-\pi }^{\pi }{\sum\nolimits_{k=1}^{\infty }{{{h}_{k}}({{e}^{i\lambda }})(Z_{k}^{\zeta +\theta }(\,d\lambda )}}).
\end{equation}
The mean square error of this estimate $\hat{A}_{N}^{i}\zeta $ is calculated by the formula
$$\Delta (h;f,g)=E\,|A_{N}^{i}\zeta -\hat{A}_{N}^{i}\zeta {{|}^{2}}=$$
\begin{equation} \label{eq08}
=\frac{1}{2\pi }\int_{-\pi }^{\pi }{\left( {{[{{A}_{N}}({{e}^{i\lambda }})-h({{e}^{i\lambda }})]}^{\top }}f(\lambda )\overline{[{{A}_{N}}({{e}^{i\lambda }})-h({{e}^{i\lambda }})]}+{{h}^{\top}}({{e}^{i\lambda }})g(\lambda )\overline{h({{e}^{i\lambda }})} \right)}\,d\lambda,
\end{equation}
$${{A}_{N}}({{e}^{i\lambda }})=\sum\nolimits_{j=0}^{N}{{{{\vec{a}}}_{j}}{{e}^{ij\lambda }}}.$$
The spectral characteristic $h(f,g)$ of the optimal linear estimate $\hat{A}_{N}^{i}\zeta $ for given spectral densities $f(\lambda ),\,\,g(\lambda )$ minimizes the value of the mean square error
\begin{equation} \label{eq09}
\Delta (f,g)=\Delta (h(f,g);f,g)=\underset{h\in L_{2}^{\,N-}(f+g)}{\mathop{\min }}\,\Delta (h;f,g)=\underset{\hat{A}_{N}^{i}\zeta }{\mathop{\min }}\,E\,|A_{N}^{i}\zeta -\hat{A}_{N}^{i}\zeta {{|}^{2}}. \end{equation}
The optimal linear estimate $\hat{A}_{N}^{i}\zeta $ is a solution of the optimization problem (9). To find the spectral characteristic $h(f,g)$ and the mean square error $\Delta (f,g)$ of the optimal linear estimate $\hat{A}_{N}^{i}\zeta $ we use the Hilbert space orthogonal projection method proposed by Kolmogorov (1992). According to the method $\hat{A}_{N}^{i}\zeta $ is a projection of $A_{N}^{i}\zeta $ on the subspace ${{H}^{N-}}(\zeta +\theta )$ generated in the space $H(\zeta +\theta )$ by values ${{\zeta }_{j}}+{{\theta }_{j}},\,j\in \mathbb{Z}\backslash \{0,1,...,N\}.$ The optimal estimate $\hat{A}_{N}^{i}\zeta $ is determined by two conditions:

1) $\hat{A}_{N}^{i}\zeta \in {{H}^{N-}}(\zeta +\theta ),$

2) $(A_{N}^{i}\zeta -\hat{A}_{N}^{i}\zeta ) \perp {{H}^{N-}}(\zeta +\theta ).$

The second condition gives the formula for the spectral characteristic $h(f,g)$ of the optimal estimate $\hat{A}_{N}^{i}\zeta $
\begin{equation}\label{eq10}
\begin{split}
{{h}^{\top}}(f,g)=\left( A_{N}^{\top }({{e}^{i\lambda }})f(\lambda )-C_{N}^{\top}({{e}^{i\lambda }}) \right){{[f(\lambda )+g(\lambda )]}^{-1}}=\\
=A_{N}^{\top }({{e}^{i\lambda }})-\left( A_{N}^{\top}({{e}^{i\lambda }})g(\lambda )+C_{N}^{\top}({{e}^{i\lambda }}) \right){{[f(\lambda )+g(\lambda )]}^{-1}},
\end{split}
\end{equation}
$${{C}_{N}}({{e}^{i\lambda }})=\sum\nolimits_{j=0}^{N}{{{{\vec{c}}}_{j}}{{e}^{ij\lambda }}}.$$

The first condition leads to the equation for unknown coefficients ${{c}_{N}}=\{{{\vec{c}}_{j}}\}_{j=0}^{N}\,$
\begin{equation}\label{eq11}
{{c}_{N}}={{B}_{N}}^{-1}{{D}_{N}}{{a}_{N}},
\end{equation}
where ${{a}_{N}}=\{{{\vec{a}}_{j}}\}_{j=0}^{N}\,$ is a vector, block-matrices ${{B}_{N}}=\{{{B}_{N}}(l,j)\}_{l,j=0}^{N}\,,$ ${{D}_{N}}=\{{{D}_{N}}(l,j)\}_{l,j=0}^{N}\,$ are determined by elements
$${{B}_{N}}(l,j)=\frac{1}{2\pi }\int_{-\pi }^{\pi }{{{\left[ {{(f(\lambda )+g(\lambda ))}^{-1}} \right]}^{\top }}\,{{e}^{i(j-l)\lambda }}d\lambda \,,}$$
$${{D}_{N}}(l,j)=\frac{1}{2\pi }\int_{-\pi }^{\pi }{{{\left[ f(\lambda ){{(f(\lambda )+g(\lambda ))}^{-1}} \right]}^{\top }}\,{{e}^{i(j-l)\lambda }}d\lambda \,.\,}$$

Taking into account formula (8) and the derived relations (10), (11), the mean square error of the optimal estimate $\hat{A}_{N}^{i}\zeta $ can be calculated by the formula
\begin{equation}\label{eq12}
\Delta (f,g)=\Delta (h(f,g);f,g)=\left\langle {{a}_{N}},{{R}_{N}}{{a}_{N}} \right\rangle +\left\langle {{c}_{N}},{{B}_{N}}{{c}_{N}} \right\rangle,
\end{equation}
where $\left\langle \cdot \,,\cdot  \right\rangle $ is the scalar product in ${{\ell }_{2}}$, the block-matrix ${{R}_{N}}=\{{{R}_{N}}(l,j)\}_{l,j=0}^{N}\,$ is determined by elements
$${{R}_{N}}(l,j)=\frac{1}{2\pi }\int_{-\pi }^{\pi }{{{\left[ f(\lambda ){{(f(\lambda )+g(\lambda ))}^{-1}}g(\lambda ) \right]}^{\top }}\,{{e}^{i(j-l)\lambda }}d\lambda \,.\,}$$
Thus our results can be summarized in the following statements. For more details see article by Dubovetska and Moklyachuk (2012c).

\begin{thm} \label{thm3.1.1}
Let $\{\zeta (t),\,t\in \mathbb{R}\}$ and $\{\theta (t),\,t\in \mathbb{R}\}$ be uncorrelated PC stochastic processes such that the generated stationary sequences $\{{{\zeta }_{j}},\,j\in \mathbb{Z}\}$ and $\{{{\theta }_{j}},\,j\in \mathbb{Z}\}$ have spectral densities $f(\lambda )$ and $g(\lambda )$, respectively, which satisfy the minimality condition (6).
Let coefficients $\{{{\vec{a}}_{j}},\,j=0,1,...,N\}$, that determine the functional $A_{N}^{i}\zeta $,  satisfy condition (5).
The spectral characteristic $h(f,g)$ and the mean square error $\Delta (f,g)$ of the optimal linear estimate of the functional $A_{N}^{i}\zeta $ based on observations of the process $\zeta (t)+\theta (t)$ at points of time $t\in \mathbb{R}\backslash [0,(N+1)T]$ are calculated by formulas (10) and (12). The optimal linear estimate $\hat{A}_{N}^{i}\zeta $is determined by the formula (7).
\end{thm}

For the problem of optimal linear estimation of the functional $A_{N}^{i}\zeta $ based on observations of the process $\zeta (t)$ without noise we have the following corollary from Theorem 3.1.1.

\begin{nas} \label{corr3.1.1}
Let $\{\zeta (t),\,t\in \mathbb{R}\}$ be PC stochastic process such that the generated stationary sequence $\{{{\zeta }_{j}},\,j\in \mathbb{Z}\}$ has the spectral density $f(\lambda )$ , which satisfy the minimality condition
\begin{equation}\label{eq13}
\int_{-\pi }^{\pi }{Tr[{{(f(\lambda ))}^{-1}}]\,d\lambda <\infty \,.}
\end{equation}
Let coefficients $\{{{\vec{a}}_{j}},\,j=0,1,...,N\}$, that determine the functional $A_{N}^{i}\zeta $,   satisfy condition (5).
The spectral characteristic $h(f)$ and the mean square error $\Delta (f)$ of the optimal linear estimate of the functional $A_{N}^{i}\zeta $ based on observations of the process $\zeta (t)$ at points of time $t\in \mathbb{R}\backslash [0,(N+1)T]$ are calculated by formulas
\begin{equation}\label{eq14}
{{h}^{\top }}(f)=A_{N}^{\top }({{e}^{i\lambda }})-C_{N}^{\top }({{e}^{i\lambda }}){{[f(\lambda )]}^{-1}}\,,
\end{equation}
\begin{equation}\label{eq15}
\Delta (f)=\left\langle {{c}_{N}},{{a}_{N}} \right\rangle,
\end{equation}
where ${{c}_{N}}=\{{{\vec{c}}_{j}}\}_{j=0}^{N}\,$, ${{c}_{N}}=B_{N}^{-1}{{a}_{N}},$ block-matrix ${{B}_{N}}=\{{{B}_{N}}(l,j)\}_{l,j=0}^{N}\,$ is determined by elements
$${{B}_{N}}(l,j)=\frac{1}{2\pi }\int_{-\pi }^{\pi }{{{\left[ {{(f(\lambda ))}^{-1}} \right]}^{\top }}\,{{e}^{i(j-l)\lambda }}d\lambda \,}.$$
 The optimal linear estimate $\hat{A}_{N}^{i}\zeta $ of the functional $A_{N}^{i}\zeta $ is determined by the formula
 \begin{equation}\label{eq16}
\hat{A}_{N}^{i}\zeta =\int_{-\pi }^{\pi }{{{h}^{\top }}({{e}^{i\lambda }})}({{Z}^{\zeta }}(\,d\lambda ))=\int_{-\pi }^{\pi }{\sum\nolimits_{k=1}^{\infty }{{{h}_{k}}({{e}^{i\lambda }})(Z_{k}^{\zeta }(\,d\lambda )}}).
\end{equation}
\end{nas}

\subsection{Extrapolation problem}

Consider the problem of optimal linear estimation of the functional
$${{A}^{e}}\zeta =\int_{0}^{\infty }{a(t)\zeta (t)\,dt}$$
that depends on the unknown values of PC stochastic process $\zeta (t)$ based on observations of the process $\zeta (t)+\theta (t)$ at points of time $t<0$. The noise $\theta (t)$ is an uncorrelated with $\zeta (t)$ PC stochastic process. The function $a(t),\,\,t\in {{\mathbb{R}}_{+}},$ satisfies condition
$$\int_{0}^{\infty }{|a(t)|\,dt}<\infty .$$
Taking into consideration transformation (1) of the process $\zeta (t)$ and decomposition (3) of the generated stationary sequence $\{{{\zeta }_{j}},\,j\in \mathbb{Z}\},$ the functional ${{A}^{e}}\zeta $can be represented in the form
$${{A}^{e}}\zeta =\int_{0}^{\infty }{a(t)\zeta (t)\,dt=}\sum\nolimits_{j=0}^{\infty }{\int_{0}^{T}{{{a}_{j}}(u){{\zeta }_{j}}(u)}\,du}=\sum\nolimits_{j=0}^{\infty }{\sum\nolimits_{k=1}^{\infty }{{{a}_{kj}}{{\zeta }_{kj}}}=\sum\nolimits_{j=0}^{\infty }{\vec{a}_{j}^{\top }{{{\vec{\zeta }}}_{j}}}},$$
where vector $\,\,{{\vec{a}}_{j}}={{({{a}_{kj}},\,k=1,2,...)}^{\top }}={{({{a}_{1j}},{{a}_{3j,}}{{a}_{2j}},...,{{a}_{2k+1,j}},{{a}_{2k,j}},...)}^{\top }}.$

Assume that coefficients $\{{{\vec{a}}_{j}},\,j=0,1,...\}$ satisfy conditions
\begin{equation}\label{eq17}
\sum\nolimits_{j=0}^{\infty }{||{{{\vec{a}}}_{j}}||}<\infty,\quad  {{\sum\nolimits_{j=0}^{\infty }{(j+1)||{{{\vec{a}}}_{j}}||}}^{2}}<\infty.
\end{equation}
It follows from the first condition from (17) that the functional ${{A}^{e}}\zeta $ has finite second moment. The second condition provides compactness of operators defined further.

Let spectral densities $f(\lambda )$ and $g(\lambda )$ of the generated stationary sequences $\{{{\zeta }_{j}},\,j\in \mathbb{Z}\}$ and $\{{{\theta }_{j}},\,j\in \mathbb{Z}\}$ be such that the minimality condition (6) is satisfied.

Denote by $L_{2}^{\,-}(f)$ the subspace of ${{L}_{2}}(f)$ generated by vector-functions ${{e}^{ij\lambda }}{{\delta }_{k}},$ $j<0,$ ${{\delta }_{k}}=\{{{\delta }_{kn}}\}_{n=1}^{\infty },$ $k=1,2,...,$ where ${{\delta }_{kn}}$ is the Kronecker symbol: ${{\delta }_{kk}}=1$ and ${{\delta }_{kn}}=0$ for $k\ne n$.
Every linear estimate ${{\hat{A}}^{e}}\zeta $ of the functional ${{A}^{e}}\zeta $ based on observations of the process $\zeta (t)+\theta (t)$ at points of time $t<0$ is determined by spectral characteristic $h({{e}^{i\lambda }})\in L_{2}^{\,-}(f+g)$ and by formula (7).

The classical Kolmogorov projection method (1992) allows us to find the value of the mean square error $\Delta (f,g)$ and spectral characteristic $h(f,g)$ of the optimal linear estimate of the functional ${{A}^{e}}\zeta $ under the condition that spectral densities $f(\lambda ),\,\,g(\lambda )$ of generated stationary sequences $\{{{\zeta }_{j}},\,j\in \mathbb{Z}\}$, $\{{{\theta }_{j}},\,j\in \mathbb{Z}\}$ are known. Applying the same considerations as in the case of interpolation problem we can verify validity of the following statements. For more details see article by Dubovetska and Moklyachuk (2013b).

\begin{thm} \label{thm3.2.1}
Let $\{\zeta (t),\,t\in \mathbb{R}\}$ and $\{\theta (t),\,t\in \mathbb{R}\}$ be uncorrelated PC stochastic processes such that the generated stationary sequences $\{{{\zeta }_{j}},\,j\in \mathbb{Z}\}$ and $\{{{\theta }_{j}},\,j\in \mathbb{Z}\}$ have spectral densities $f(\lambda )$ and $g(\lambda )$, respectively, which satisfy the minimality condition (6). Let coefficients $\{{{\vec{a}}_{j}},\,j=0,1,...\}$, that determine the functional ${{A}^{e}}\zeta $,  satisfy conditions (17).
The spectral characteristic $h(f,g)$ and the mean square error $\Delta (f,g)$ of the optimal linear estimate of the functional ${{A}^{e}}\zeta $ based on observations of the process $\zeta (t)+\theta (t)$ at points of time $t<0$ are calculated by formulas

$${{h}^{\top }}(f,g)=\left( {{A}^{\top }}({{e}^{i\lambda }})f(\lambda )-{{C}^{\top }}({{e}^{i\lambda }}) \right){{[f(\lambda )+g(\lambda )]}^{-1}}=$$
\begin{equation}\label{eq18}
={{A}^{\top }}({{e}^{i\lambda }})-\left( {{A}^{\top }}({{e}^{i\lambda }})g(\lambda )+{{C}^{\top }}({{e}^{i\lambda }}) \right){{[f(\lambda )+g(\lambda )]}^{-1}},
\end{equation}
\begin{equation}\label{eq19}
\Delta (f,g)=\Delta (h(f,g);f,g)=\left\langle a,Ra \right\rangle +\left\langle c,Bc \right\rangle,
\end{equation}
where $A({{e}^{i\lambda }})=\sum\nolimits_{j=0}^{\infty }{{{{\vec{a}}}_{j}}{{e}^{ij\lambda }}},$ $C({{e}^{i\lambda }})=\sum\nolimits_{j=0}^{\infty }{{{{\vec{c}}}_{j}}{{e}^{ij\lambda }}},$ vector $a=\{{{\vec{a}}_{j}}\}_{j=0}^{\infty }\,$ is given, vector of unknown coefficients $c=\{{{\vec{c}}_{j}}\}_{j=0}^{\infty }\,$ is determined by the equation $c={{B}^{-1}}Da,$ block-matrices $B=\{B(l,j)\}_{l,j=0}^{\infty }\,,$ $D=\{D(l,j)\}_{l,j=0}^{\infty }\,,$$R=\{R(l,j)\}_{l,j=0}^{\infty }\,$ are determined by elements
$$B(l,j)=\frac{1}{2\pi }\int_{-\pi }^{\pi }{{{\left[ {{(f(\lambda )+g(\lambda ))}^{-1}} \right]}^{\top }}\,{{e}^{i(j-l)\lambda }}d\lambda \,,}$$
$$D(l,j)=\frac{1}{2\pi }\int_{-\pi }^{\pi }{{{\left[ f(\lambda ){{(f(\lambda )+g(\lambda ))}^{-1}} \right]}^{\top }}\,{{e}^{i(j-l)\lambda }}d\lambda \,,\,}$$
$$R(l,j)=\frac{1}{2\pi }\int_{-\pi }^{\pi }{{{\left[ f(\lambda ){{(f(\lambda )+g(\lambda ))}^{-1}}g(\lambda ) \right]}^{\top }}\,{{e}^{i(j-l)\lambda }}d\lambda \,,\,\,\,l,j=0,1,...\,}$$
The optimal linear estimate ${{\hat{A}}^{e}}\zeta $ of the functional ${{A}^{e}}\zeta $ is determined by the formula (7).
\end{thm}

For the problem of optimal linear estimation of the functional ${{A}^{e}}\zeta $ based on observations of the process $\zeta (t)$ without noise we have the following statement, which is a corollary from Theorem \ref{thm3.2.1}.

\begin{nas} \label{corr3.2.1}
  Let $\{\zeta (t),\,t\in \mathbb{R}\}$ be PC stochastic process such that the generated stationary sequence $\{{{\zeta }_{j}},\,j\in \mathbb{Z}\}$ has spectral density $f(\lambda )$, which satisfy the minimality condition (13). Let  coefficients $\{{{\vec{a}}_{j}},\,j=0,1,...\}$, that determine the functional ${{A}^{e}}\zeta $,   satisfy conditions (17).
The spectral characteristic $h(f)$ and the mean square error $\Delta (f)$ of the optimal linear estimate ${{\hat{A}}^{e}}\zeta $ of the functional ${{A}^{e}}\zeta $ based on observations of the process $\zeta (t)$ at points of time $t<0$ are calculated by formulas
\begin{equation}\label{eq20}
{{h}^{\top }}(f)={{A}^{\top }}({{e}^{i\lambda }})-{{C}^{\top }}({{e}^{i\lambda }}){{[f(\lambda )]}^{-1}},
\end{equation}
\begin{equation}\label{eq21}
\Delta (f)=\left\langle c,a \right\rangle,
\end{equation}
where $c=\{{{\vec{c}}_{j}}\}_{j=0}^{\infty }\,$, $c={{B}^{-1}}a,$ block-matrix $B=\{B(l,j)\}_{l,j=0}^{\infty }\,$ is determined by elements
$$B(l,j)=\frac{1}{2\pi }\int_{-\pi }^{\pi }{{{\left[ {{(f(\lambda ))}^{-1}} \right]}^{\top }}\,{{e}^{i(j-l)\lambda }}d\lambda \,}.$$
The optimal linear estimate ${{\hat{A}}^{e}}\zeta $ of the functional ${{A}^{e}}\zeta $ is determined by the formula  (16).
\end{nas}

Note that Kolmogorov proposed a method of solving the problem of interpolation of stationary sequence (i. e. finding spectral characteristic and mean square error of the optimal linear estimate of one missed observation of the sequence) using the Fourier coefficients of the function ${1}/{f}\;$. Theorem 3.2.1 shows that the Fourier coefficients of some functions from spectral densities can be used to find spectral characteristics and the mean square error of optimal linear estimates of functionals of  stationary sequences for problems of extrapolation and interpolation based on observations without noise as well as on observations with noise.

The form of the spectral characteristics and the form of the mean square error of the optimal linear estimate are convenient for finding the least favourable spectral densities and minimax spectral characteristics of optimal estimates for the problems of extrapolation and interpolation based on observations without noise as well as on observations with noise.

To solve the problem of extrapolation of stationary sequences Kolmogorov (see also Kailath (1974), Rozanov (1967), Wiener (1966), Yaglom (1987)) proposed a method based on factorization of spectral density. This method is suitable for solving problems of extrapolation based on observations without noise whereas Theorem 3.2.1 describes the method of solving problem of extrapolation based on observations with noise.

Let apply the method based on factorization of the spectral density to the problem of estimation of the functional from observations without noise. For more results see articles by Moklyachuk (1995,1996) and book by Moklyachuk (2008).

\begin{ozn} \label{def3.2.1}
 Denote by ${{H}_{\zeta }}(n)$ the closed linear subspace of the Hilbert space $H={{L}_{2}}(\Omega ,F,P)$ generated by random variables $\{{{\zeta }_{kj}},\,k\ge 1\,,\,j\le n\}.$ The sequence $\{{{\zeta }_{j}},\,j\in \mathbb{Z}\}$ is called regular if $\bigcap\nolimits_{n}{{{H}_{\zeta }}(n)=}\varnothing .$ If $\bigcap\nolimits_{n}{{{H}_{\zeta }}(n)=}H$ then the sequence $\{{{\zeta }_{j}},\,j\in \mathbb{Z}\}$ is called singular.
\end{ozn}

The following result allows simplifying the problem of optimal linear estimation of unknown values of stationary sequence.

\begin{thm} \label{thm3.2.2}
 A stationary sequence $\{{{\zeta }_{j}},\,j\in \mathbb{Z}\}$ admits a unique representation in the form
${{\zeta }_{j}}=\zeta _{j}^{r}+\zeta _{j}^{s}$
where $\{\zeta _{j}^{r},\,j\in \mathbb{Z}\}$ is a regular sequence and $\{\zeta _{j}^{s},\,j\in \mathbb{Z}\}$ is a singular sequence.
Moreover, the sequences $\{\zeta _{j}^{r}\}$ and $\{\zeta _{n}^{s}\}$ are orthogonal for all $j,n\in \mathbb{Z}$.
\end{thm}

Since the unknown values of components of singular stationary sequence has error-free estimate, we will consider the estimation problem only for regular stationary sequences.

The regular stationary sequence $\{{{\zeta }_{j}},\,j\in \mathbb{Z}\}$ admits the canonical moving average representation of components (Kallianpur and Mandrekar, 1971; Moklyachuk, 1981)
\begin{equation}\label{eq22}
{{\zeta }_{kj}}=\sum\nolimits_{u=-\infty }^{j}{\sum\nolimits_{m=1}^{M}{{{d}_{km}}(j-u){{\varepsilon }_{m}}(u)}},
\end{equation}
where ${{\varepsilon }_{m}}(u),\,m=1,...,M,\,u\in \mathbb{Z}$ are mutually orthogonal sequences in $H$ with orthonormal values; $M$ is the multiplicity of $\{{{\zeta }_{j}},\,j\in \mathbb{Z}\};$ ${{d}_{km}}(u),\,\,k=1,2,...,\,\,\,m=1,...,\,\,M,\,\,\,u=0,1,...,$ are matrix-valued sequences such that $\sum\nolimits_{u=0}^{\infty }{\sum\nolimits_{k=1}^{\infty }{\sum\nolimits_{m=1}^{M}{|{{d}_{km}}(u){{|}^{2}}}}}={{P}_{\zeta }}.$

As a consequence of representation (22) the optimal linear estimate of components of stationary sequence $\{{{\zeta }_{j}},\,j\in \mathbb{Z}\}$ can be represented in the form
\begin{equation}\label{eq23}
{{\hat{\zeta }}_{kj}}=\sum\nolimits_{s=-\infty }^{-1}{\sum\nolimits_{m=1}^{M}{{{d}_{km}}(j-s){{\varepsilon }_{m}}(s)}}.
\end{equation}
The spectral density $f(\lambda )$ of regular stationary sequence $\{{{\zeta }_{j}},\,j\in \mathbb{Z}\}$ admits the canonical factorization
\begin{equation}\label{eq22}
f(\lambda )=P(\lambda )\,{{P}^{*}}(\lambda ),\quad P(\lambda )=\sum\nolimits_{u=0}^{\infty }{d(u){{e}^{-iu\lambda }}},
\end{equation}
where matrices $d(u)=\{{{d}_{km}}(u)\}_{k=\overline{1,\infty }}^{m=\overline{1,M}}$, $u\ge 0,$ are determined by coefficients of the canonical representation (22).

Taking into account decompositions (22), (23) and factorization (24) we can verify validity of the following result.

\begin{thm} \label{thm3.2.3}
Let $\{\zeta (t),\,t\in \mathbb{R}\}$ be a PC stochastic process such that the generated stationary sequence $\{{{\zeta }_{j}},\,j\in \mathbb{Z}\}$ has  spectral density $f(\lambda )$ , which satisfy the minimality condition (13) and admits the canonical factorization (24). Let coefficients $\{{{\vec{a}}_{j}},\,j=0,1,...\}$, that determine the functional ${{A}^{e}}\zeta $,   satisfy conditions (17).
The spectral characteristic $h(f)$ and the mean square error $\Delta (f)$ of the optimal linear estimate of the functional ${{A}^{e}}\zeta $ based on observations of the process $\zeta (t)$ at points of time $t<0$ are calculated by formulas
\begin{equation}\label{eq25}
{{h}^{\top }}(f)={{A}^{\top }}({{e}^{i\lambda }})-S({{e}^{i\lambda }})Q(\lambda )\,,
\end{equation}
\begin{equation}\label{eq26}
\Delta (f)={{\left\| Ad \right\|}^{2}},
\end{equation}
where $S({{e}^{i\lambda }})=\sum\nolimits_{l=0}^{\infty }{{{(Ad)}_{l}}{{e}^{il\lambda }}},$ ${{(Ad)}_{l}}=\sum\nolimits_{j=l}^{\infty }{\vec{a}_{j}^{\top }d(j-l)},\,\,\,\,l\ge 0,$ the matrix function $Q(\lambda )=\{{{q}_{mk}}(\lambda )\}_{m=\overline{1,M}}^{k=\overline{1,\infty }}$ satisfies equation $Q(\lambda )P(\lambda )={{I}_{M}},$ ${{\left\| Ad \right\|}^{2}}={{\sum\nolimits_{l=0}^{\infty }{\left\| {{(Ad)}_{l}} \right\|}}^{2}}.$ The optimal linear estimate ${{\hat{A}}^{e}}\zeta $ of the functional ${{A}^{e}}\zeta $ is determined by  formula (16).
\end{thm}

Similar reasoning can be applied to find the optimal estimate of the functional
$$A_{N}^{e}\zeta =\int_{0}^{(N+1)T}{a(t)\zeta (t)\,dt=}\sum\nolimits_{j=0}^{N}{\int_{0}^{T}{{{a}_{j}}(u){{\zeta }_{j}}(u)}\,du}=\sum\nolimits_{j=0}^{N}{\sum\nolimits_{k=1}^{\infty }{{{a}_{kj}}{{\zeta }_{kj}}}=\sum\nolimits_{j=0}^{N}{\vec{a}_{j}^{\top }{{{\vec{\zeta }}}_{j}}}}.$$

The following corollary from Theorem \ref{thm3.2.3} holds true.

\begin{nas} \label{corr3.2.2}
 Let $\{\zeta (t),\,t\in \mathbb{R}\}$ be a PC stochastic process such that the generated stationary sequence $\{{{\zeta }_{j}},\,j\in \mathbb{Z}\}$ has spectral density $f(\lambda )$ , which satisfy the minimality condition (13) and admits the canonical factorization (24). Let coefficients $\{{{\vec{a}}_{j}},\,j=0,1,...,N\}$, that determine the functional $A_{N}^{e}\zeta $,   satisfy condition (5). The spectral characteristic ${{h}_{N}}(f)$ and the mean square error ${{\Delta }_{N}}(f)$ of the optimal linear estimate of the functional $A_{N}^{e}\zeta $ based on observations of the process $\zeta (t)$ at points of time $t<0$ are calculated by formulas
\begin{equation}\label{eq28}
h_{N}^{\top }(f)=A_{N}^{\top }({{e}^{i\lambda }})-{{S}_{N}}({{e}^{i\lambda }})Q(\lambda )\,,
\end{equation}
\begin{equation}\label{eq28}
{{\Delta }_{N}}(f)={{\sum\nolimits_{l=0}^{\infty }{\left\| {{({{A}_{N}}d)}_{l}} \right\|}}^{2}}={{\left\| {{A}_{N}}d \right\|}^{2}},
\end{equation}
where ${{A}_{N}}({{e}^{i\lambda }})=\sum\nolimits_{j=0}^{N}{{{{\vec{a}}}_{j}}{{e}^{ij\lambda }}},$ ${{S}_{N}}({{e}^{i\lambda }})=\sum\nolimits_{l=0}^{N}{{{({{A}_{N}}d)}_{l}}{{e}^{il\lambda }}},$ ${{({{A}_{N}}d)}_{l}}=\sum\nolimits_{j=l}^{N}{\vec{a}_{j}^{\top }d(j-l)}.$
The optimal linear estimate $\hat{A}_{N}^{e}\zeta $ of the functional $A_{N}^{e}\zeta $ is determined by the formula (16).
\end{nas}

For a specified class of PC processes we can calculate the greatest values of the mean square error $\Delta (\zeta ,\hat{A}_{{}}^{e})=E|{{A}^{e}}\zeta -{{\hat{A}}^{e}}\zeta {{|}^{2}}$ of estimate ${{\hat{A}}^{e}}\zeta $ of the functional ${{A}^{e}}\zeta $ and of the mean square error $\Delta (\zeta ,\hat{A}_{N}^{e})=E|A_{N}^{e}\zeta -\hat{A}_{N}^{e}\zeta {{|}^{2}}$ of estimate $\hat{A}_{N}^{e}\zeta $ of the functional $A_{N}^{e}\zeta $. For proof of the following results see the article by Dubovetska and Moklyachuk (2013a).

Denote by $\Lambda $ the set of all linear estimates of the functional ${{A}^{e}}\zeta $ based on observation of the process $\zeta (t)$ at points of time $t<0$. Let $\mathbf{Y}$ denotes the class of mean square continuous PC processes $\zeta (t)$ such that $E\zeta (t)=0$ and $E|\zeta (t){{|}^{2}}\le {{{P}_{\zeta }}}/{T}\;.$ The following theorem holds true.

\begin{thm} \label{thm3.2.4}

Let coefficients $\{{{\vec{a}}_{j}},\,j=0,1,...,N\}$ which determine the functional $A_{N}^{e}\zeta $ satisfy condition (5). The function $\Delta (\zeta ,\hat{A}_{N}^{e})$ has a saddle point on the set $\mathbf{Y}\times \Lambda $ and the following equality holds true
$$\underset{\hat{A}_{N}^{e}\in \Lambda }{\mathop{\min }}\,\underset{\zeta \in \mathbf{Y}}{\mathop{\max }}\,\Delta (\zeta ,\hat{A}_{N}^{e})=\underset{\zeta \in \mathbf{Y}}{\mathop{\max }}\,\underset{\hat{A}_{N}^{e}\in \Lambda }{\mathop{\min }}\,\Delta (\zeta ,\hat{A}_{N}^{e})={{P}_{\zeta }}\cdot \nu _{N}^{2},$$
where $\nu _{N}^{2}$ is the greatest eigenvalue of the self-adjoint compact operator ${{Q}_{N}}=\{{{Q}_{N}}(p,q)\}_{p,q=0}^{N}$ in the space ${{\ell }_{2}}$ determined by block-matrices ${{Q}_{N}}(p,q)=\{Q_{kn}^{N}(p,q)\}_{k,n=1}^{\infty }$ with elements
$$Q_{kn}^{N}(p,q)=\sum\nolimits_{s=0}^{min(N-p,N-q)}{{{a}_{k,s+p}}\cdot \overline{{{a}_{n,s+q}}}}, $$
$$\,\,\,k,n=1,2,\ldots ,\quad p,q=0,1,\ldots ,N.$$
The least favorable stochastic sequence generated by PC process from the class $\mathbf{Y}$ for the optimal estimate of the functional $A_{N}^{e}\zeta $ is a one-sided moving average sequence of order $N$ of the form
$${{\vec{\zeta }}_{j}}=\sum\nolimits_{u=j-N}^{j}{d(j-u)\vec{\varepsilon }(u)},$$
where ${{d}_{N}}=(d(p))_{p=0}^{N}$ is the eigenvector, that corresponds to $\nu _{N}^{2}$, which is constructed from matrices $d(p)=\{{{d}_{km}}(p)\}_{k=\overline{1,\infty }}^{m=\overline{1,M}}$ and is determined by condition $||{{d}_{N}}|{{|}^{2}}=\sum\nolimits_{p=0}^{N}{||d(p)|{{|}^{2}}}={{P}_{\zeta }}$, $\vec{\varepsilon }(u)=\{{{\varepsilon }_{m}}(u)\}_{m=1}^{M}$  is a vector-valued stationary stochastic sequence with orthogonal values.
\end{thm}

\begin{thm} \label{thm3.2.5}
 Let coefficients $\{{{\vec{a}}_{j}},\,j=0,1,...\}$ which determine the functional ${{A}^{e}}\zeta $ satisfy conditions (17). The function $\Delta (\zeta ,\hat{A}_{{}}^{e})$ has a saddle point on the set $\mathbf{Y}\times \Lambda $ and the following equality holds true
$$\underset{{{{\hat{A}}}^{e}}\in \Lambda }{\mathop{\min }}\,\underset{\zeta \in \mathbf{Y}}{\mathop{\max }}\,\Delta (\zeta ,{{\hat{A}}^{e}})=\underset{\zeta \in \mathbf{Y}}{\mathop{\max }}\,\underset{{{{\hat{A}}}^{e}}\in \Lambda }{\mathop{\min }}\,\Delta (\zeta ,{{\hat{A}}^{e}})={{P}_{\zeta }}\cdot {{\nu }^{2}},$$
where $\nu _{{}}^{2}$ is the greatest eigenvalue of the self-adjoint compact operator $Q=\{Q(p,q)\}_{p,q=0}^{\infty }$ in the space ${{\ell }_{2}}$ determined by block-matrices $Q(p,q)=\{{{Q}_{kn}}(p,q)\}_{k,n=1}^{\infty }$ with elements
$$Q_{kn}^{{}}(p,q)=\sum\nolimits_{s=0}^{\infty }{{{a}_{k,s+p}}\cdot \overline{{{a}_{n,s+q}}}},$$
$$k,n=1,2,\ldots ,\quad p,q=0,1,\ldots .$$
The least favorable stochastic sequence generated by PC process from the class $\mathbf{Y}$ for the optimal estimate of the functional ${{A}^{e}}\zeta $ is an one-sided moving average sequence of the form
$${{\vec{\zeta }}_{j}}=\sum\nolimits_{u=-\infty }^{j}{d(j-u)\vec{\varepsilon }(u)},$$
where $d=(d(p))_{p=0}^{\infty }$ is the eigenvector, that corresponds to $\nu _{{}}^{2}$, it is constructed from matrices $d(p)=\{{{d}_{km}}(p)\}_{k=\overline{1,\infty }}^{m=\overline{1,M}}$ and is determined by the condition  $||d|{{|}^{2}}=\sum\nolimits_{p=0}^{\infty }{||d(p)|{{|}^{2}}}={{P}_{\zeta }}.$
\end{thm}

\subsection{Filtering problem}

Consider the problem of optimal linear estimation of the functional
$$A_{{}}^{f}\zeta =\int_{0}^{\infty }{a}(t)\zeta (-t)dt$$
which depends on the unknown values of a PC  stochastic process $\zeta (t)$ based on observations of the process $\zeta (t)+\theta (t)$ at points of time $t\le 0$. The function $a(t),\,\,t\in {{\mathbb{R}}_{+}},$ satisfies the condition $\int_{0}^{\infty }{|a(t)|\,dt}<\infty .$

With the help of transformations (1), (3) of the process $\zeta (t)$ we can represent the functional $A_{{}}^{f}\zeta $ in the form
$${{A}^{f}}\zeta =\int_{0}^{\infty }{a(t)\zeta (-t)\,dt=}\sum\nolimits_{j=0}^{\infty }{\int_{0}^{T}{{{a}_{j}}(u){{\zeta }_{-j}}(-u)}\,du}=\sum\nolimits_{j=0}^{\infty }{\sum\nolimits_{k=1}^{\infty }{{{a}_{kj}}{{\zeta }_{k,-j}}}=\sum\nolimits_{j=0}^{\infty }{\vec{a}_{j}^{\top }{{{\vec{\zeta }}}_{-j}}}},$$
where ${{a}_{j}}(u)=a(u+jT),$  ${{\zeta }_{-j}}(-u)=\zeta (-u-jT),$  $u\in [0,T),$ and vector ${{\vec{a}}_{j}}={{({{a}_{kj}},\,k=1,2,...)}^{\top }}={{({{a}_{1j}},{{a}_{3j,}}{{a}_{2j}},...,{{a}_{2k+1,j}},{{a}_{2k,j}},...)}^{\top }}.$

Assume that coefficients $\{{{\vec{a}}_{j}},\,j=0,1,...\}$ satisfy condition
\begin{equation}\label{eq29}
\sum\nolimits_{j=0}^{\infty }{||{{{\vec{a}}}_{j}}||}<\infty .
\end{equation}
Condition (29) guarantees finite second moment of the functional $A_{{}}^{f}\zeta $.
Let spectral densities $f(\lambda )$ and $g(\lambda )$ of generated stationary sequences $\{{{\zeta }_{j}},\,j\in \mathbb{Z}\}$ and $\{{{\theta }_{j}},\,j\in \mathbb{Z}\}$ be such that minimality condition (6) is satisfied.

Every linear estimate $\hat{A}_{{}}^{f}\zeta $ of the functional ${{A}^{f}}\zeta $ is determined by the spectral characteristic $h({{e}^{i\lambda }})\in L_{2}^{\,-\,\,0}(f+g)$ and by formula (7). Denote by $L_{2}^{\,-\,\,0}(f+g)$ the subspace of the space $L_{2}^{{}}(f+g)$ generated by vector-functions ${{e}^{ij\lambda }}{{\delta }_{k}},$ $j\le 0,$ ${{\delta }_{k}}=\{{{\delta }_{kn}}\}_{n=1}^{\infty },$ $k=1,2,...\,\,\,.$

The mean square error of the estimate ${{\hat{A}}^{f}}\zeta $ is calculated by the formula
$$\Delta (h;f,g)=E\,|{{A}^{f}}\zeta -{{\hat{A}}^{f}}\zeta {{|}^{2}}=$$
\begin{equation}\label{eq30}
=\frac{1}{2\pi }\int_{-\pi }^{\pi }{\left( {{[{{A}_{-}}({{e}^{i\lambda }})-h({{e}^{i\lambda }})]}^{\top }}f(\lambda )\overline{[{{A}_{-}}({{e}^{i\lambda }})-h({{e}^{i\lambda }})]}+{{h}^{\top }}({{e}^{i\lambda }})g(\lambda )\overline{h({{e}^{i\lambda }})} \right)}\,d\lambda ,
\end{equation}
$${{A}_{-}}({{e}^{i\lambda }})=\sum\nolimits_{j=0}^{\infty }{{{{\vec{a}}}_{j}}{{e}^{-ij\lambda }}}.$$

The spectral characteristic $h(f,g)$ of the optimal linear estimate ${{\hat{A}}^{f}}\zeta $ for the given densities  $f(\lambda ),\,\,g(\lambda )$ minimizes the value of the mean square error
\begin{equation}\label{eq31}
\Delta (f,g)=\Delta (h(f,g);f,g)=\underset{h\in L_{2}^{\,-\,0}(f+g)}{\mathop{\min }}\,\Delta (h;f,g)=\underset{{{{\hat{A}}}^{f}}\zeta }{\mathop{\min }}\,E\,|{{A}^{f}}\zeta -{{\hat{A}}^{f}}\zeta {{|}^{2}}.
\end{equation}

The optimal linear estimate ${{\hat{A}}^{f}}\zeta $ is a solution of the optimization problem (31).
The classical Kolmogorov projection method (1992) allows us to find the spectral characteristic $h(f,g)$ and the value of the mean square error $\Delta (f,g)$ of the optimal linear estimate of the functional ${{A}^{f}}\zeta $ on condition that the spectral densities  $f(\lambda )$ and $g(\lambda )$ are known. In this case
$${{h}^{\top }}(f,g)=\left( A_{-}^{\top }({{e}^{i\lambda }})f(\lambda )-{{D}^{\top }}({{e}^{i\lambda }}) \right){{[f(\lambda )+g(\lambda )]}^{-1}}=$$
\begin{equation}\label{eq32}
=A_{-}^{\top }({{e}^{i\lambda }})-\left( A_{-}^{\top }({{e}^{i\lambda }})g(\lambda )+{{D}^{\top }}({{e}^{i\lambda }}) \right){{[f(\lambda )+g(\lambda )]}^{-1}},
\end{equation}
$$D({{e}^{i\lambda }})=\sum\nolimits_{j=1}^{\infty }{{{{\vec{d}}}_{j}}{{e}^{ij\lambda }}},$$
\begin{equation}\label{eq33}
\Delta (f,g)=\Delta (h(f,g);f,g)=\left\langle a,Wa \right\rangle +\left\langle d,Ud \right\rangle ,
\end{equation}
where vector $a=\{{{\vec{a}}_{j}}\}_{j=0}^{\infty }\,$, unknown coefficients $d=\{{{\vec{d}}_{j}}\}_{j=1}^{\infty }\,={{U}^{-1}}Va,$ block-matrices $U=\{U(j,l)\}_{l,j=1}^{\infty }\,,$ $V=\{V(j,\tilde{l})\}_{j=\overline{1,\infty }}^{\tilde{l}=\overline{0,\infty }}\,,$ $W=\{W(\tilde{j},\tilde{l})\}_{\tilde{l},\tilde{j}=0}^{\infty }\,$ are determined by elements
$$U(j,l)=\frac{1}{2\pi }\int_{-\pi }^{\pi }{{{\left[ {{(f(\lambda )+g(\lambda ))}^{-1}} \right]}^{\top }}\,{{e}^{i(l-j)\lambda }}d\lambda \,,\,\,}\,\,l,j,=1,2,...,$$
$$V(j,\tilde{l})=\frac{1}{2\pi }\int_{-\pi }^{\pi }{{{\left[ f(\lambda ){{(f(\lambda )+g(\lambda ))}^{-1}} \right]}^{\top }}\,{{e}^{-i(\tilde{l}+j)\lambda }}d\lambda \,,\,\,\,\,\,j=1,2,...,\,\,\tilde{l}=0,1,...,}$$
$$W(\tilde{j},\tilde{l})=\frac{1}{2\pi }\int_{-\pi }^{\pi }{{{\left[ f(\lambda ){{(f(\lambda )+g(\lambda ))}^{-1}}g(\lambda ) \right]}^{\top }}\,{{e}^{i(\tilde{l}-\tilde{j})\lambda }}d\lambda \,,\,\,\,\,\tilde{l},\tilde{j}=0,1,....\,}$$
Thus our results can be summarized in the following statement.

\begin{thm} \label{thm3.3.1}
Let $\{\zeta (t),\,t\in \mathbb{R}\}$ and $\{\theta (t),\,t\in \mathbb{R}\}$ be uncorrelated PC stochastic processes such that the generated stationary sequences $\{{{\zeta }_{j}},\,j\in \mathbb{Z}\}$ and $\{{{\theta }_{j}},\,j\in \mathbb{Z}\}$ have spectral densities $f(\lambda )$ and $g(\lambda )$, respectively, which satisfy the minimality condition (6). Let coefficients $\{{{\vec{a}}_{j}},\,j=0,1,...\}$, that determine the functional $A_{{}}^{f}\zeta $,  satisfy condition (29).

The spectral characteristic $h(f,g)$ and the mean square error $\Delta (f,g)$ of the optimal linear estimate of the functional $A_{{}}^{f}\zeta $ based on observations of the process $\zeta (t)+\theta (t)$ at points of time  $t\le 0$ are calculated by formulas (32) and (33), respectively. The optimal linear estimate of the functional $A_{{}}^{f}\zeta $ is determined by formula (7).
\end{thm}

\section{{Minimax-robust estimation method}}

The proposed in Section 3 formulas for calculating spectral characteristics and mean square errors of optimal linear estimates of functionals $A_{N}^{i}\zeta ,$ ${{A}^{e}}\zeta ,$ $A_{N}^{e}\zeta ,$ ${{A}^{f}}\zeta $ can be used only in the case where spectral density matrices $f(\lambda )$ and $g(\lambda )$ of the generated stationary sequences $\{{{\zeta }_{j}},\,j\in \mathbb{Z}\}$ and $\{{{\theta }_{j}},\,j\in \mathbb{Z}\}$ are exactly known. In the case where the spectral density matrices $f(\lambda )$ and $g(\lambda )$ are not exactly known, but a set $D={{D}_{f}}\times {{D}_{g}}$ of admissible spectral densities is specified, we find estimates that minimize the mean square error for all spectral densities from a given class $D$ simultaneously. Such approach to the estimation problem of functionals of the unknown values of stochastic processes is called minimax (robust) (Moklyachuk, 2008).

\begin{ozn} \label{def4.1}
 For a given class of spectral densities $D={{D}_{f}}\times {{D}_{g}}$ spectral densities ${{f}^{0}}(\lambda )\in {{D}_{f}},$ ${{g}^{0}}(\lambda )\in {{D}_{g}}$ are called least favorable in $D$ for the optimal linear estimatation of the functional $A\zeta $ if
     \[\Delta ({{f}^{0}},{{g}^{0}})=\Delta (h({{f}^{0}},{{g}^{0}});{{f}^{0}},{{g}^{0}})=\underset{(f,g)\in D}{\mathop{\max }}\,\Delta (h(f,g);f,g).\]
\end{ozn}

\begin{ozn} \label{def4.2}
 For a given class of spectral densities $D={{D}_{f}}\times {{D}_{g}}$  the spectral characteristic ${{h}^{0}}(\lambda )$ of the optimal linear estimate of the functional $A\zeta $  is called minimax-robust if
     \[{{h}^{0}}(\lambda )\in {{H}_{D}}=\bigcap\nolimits_{(f,g)\in D}{L_{2}^{\,*}(f+g)},   \underset{h\in {{H}_{D}}}{\mathop{\min }}\,\underset{(f,g)\in D}{\mathop{\max }}\,\Delta (h;f,g)=\underset{(f,g)\in D}{\mathop{\max }}\,\Delta ({{h}^{0}};f,g).\]
\end{ozn}

Here $L_{2}^{\,*}(f+g)$ denotes the subspace $L_{2}^{\,-}(f+g)$ in the case of extrapolation problem, $L_{2}^{\,N-}(f+g)$ in the case of interpolation problem, and $L_{2}^{\,-\,0}(f+g)$ in the case of  filtering problem.

Taking into consideration these definitions and Theorems 3.1.1, 3.2.1, 3.2.3, 3.3.1 we can verify that the following lemmas hold true.

\begin{lem} \label{lem4.1}
 Spectral densities ${{f}^{0}}(\lambda )\in {{D}_{f}},$ ${{g}^{0}}(\lambda )\in {{D}_{g}}$, which satisfy the minimality condition (6), are least favorable in the class $D={{D}_{f}}\times {{D}_{g}}$ for the optimal linear estimation of the functional $A_{N}^{i}\zeta $ if the Fourier coefficients of matrix functions
$${{({{f}^{0}}(\lambda )+{{g}^{0}}(\lambda ))}^{-1}},\,{{f}^{0}}(\lambda ){{({{f}^{0}}(\lambda )+{{g}^{0}}(\lambda ))}^{-1}},\,\,{{f}^{0}}(\lambda ){{({{f}^{0}}(\lambda )+{{g}^{0}}(\lambda ))}^{-1}}{{g}^{0}}(\lambda )$$
determine matrices  $B_{N}^{0},\,D_{N}^{0},\,R_{N}^{0}$, which give a solution of the extremum problem
$$\underset{(f,g)\in D}{\mathop{\max }}\,\left( \left\langle {{a}_{N}},{{R}_{N}}{{a}_{N}} \right\rangle +\left\langle {{({{B}_{N}})}^{-1}}{{D}_{N}}{{a}_{N}},{{D}_{N}}{{a}_{N}} \right\rangle  \right)=\left\langle {{a}_{N}},R_{N}^{0}{{a}_{N}} \right\rangle +\left\langle {{(B_{N}^{0})}^{-1}}D_{N}^{0}{{a}_{N}},D_{N}^{0}{{a}_{N}} \right\rangle .$$
The minimax-robust spectral characteristic ${{h}^{0}}=h({{f}^{0}},{{g}^{0}})$ of the optimal linear estimate of the functional $A_{N}^{i}\zeta $ is calculated by formula (10) if the condition $h({{f}^{0}},{{g}^{0}})\in {{H}_{D}}$ holds true.
\end{lem}
\begin{lem} \label{lem4.2}
 Spectral densities ${{f}^{0}}(\lambda )\in {{D}_{f}},$ ${{g}^{0}}(\lambda )\in {{D}_{g}}$, which satisfy the minimality condition (6), are least favorable in the class $D={{D}_{f}}\times {{D}_{g}}$ for the optimal linear estimation of the functional ${{A}^{e}}\zeta $ if the Fourier coefficients of matrix functions
$${{({{f}^{0}}(\lambda )+{{g}^{0}}(\lambda ))}^{-1}},\,\,{{f}^{0}}(\lambda ){{({{f}^{0}}(\lambda )+{{g}^{0}}(\lambda ))}^{-1}},\,\,{{f}^{0}}(\lambda ){{({{f}^{0}}(\lambda )+{{g}^{0}}(\lambda ))}^{-1}}{{g}^{0}}(\lambda )$$
determine matrices  ${{B}^{0}},\,{{D}^{0}},\,{{R}^{0}}$, which give a solution of the extremum problem
$$\underset{(f,g)\in D}{\mathop{\max }}\,\left( \left\langle a,Ra \right\rangle +\left\langle {{B}^{-1}}Da,Da \right\rangle  \right)=\left\langle a,{{R}^{0}}a \right\rangle +\left\langle {{({{B}^{0}})}^{-1}}{{D}^{0}}a,{{D}^{0}}a \right\rangle .$$
The minimax-robust spectral characteristic ${{h}^{0}}=h({{f}^{0}},{{g}^{0}})$ of the optimal linear estimate of the functional ${{A}^{e}}\zeta $ is calculated by formula (18) if the condition $h({{f}^{0}},{{g}^{0}})\in {{H}_{D}}$ holds true.
\end{lem}
\begin{lem} \label{lem4.3}
 Spectral density ${{f}^{0}}(\lambda )\in {{D}_{f}}$ is least favorable in the class ${{D}_{f}}$ for the optimal linear estimation of the functional ${{A}^{e}}\zeta $ based on observations of the process $\zeta (t)$ for $t<0$,  if it admits the canonical factorization
$${{f}^{0}}(\lambda )=\left( \sum\nolimits_{u=0}^{\infty }{{{d}^{0}}(u){{e}^{-iu\lambda }}} \right){{\left( \sum\nolimits_{u=0}^{\infty }{{{d}^{0}}(u){{e}^{-iu\lambda }}} \right)}^{*}},$$
where ${{d}^{0}}=\{{{d}^{0}}(u),\,\,u=0,1,...\}$ is a solution of the conditional extremum problem
$${{\left\| Ad \right\|}^{2}}\to \,\max,\quad f(\lambda )=\left( \sum\nolimits_{u=0}^{\infty }{d(u){{e}^{-iu\lambda }}} \right){{\left( \sum\nolimits_{u=0}^{\infty }{d(u){{e}^{-iu\lambda }}} \right)}^{*}}\in {{D}_{f}}.$$
The minimax-robust spectral characteristic ${{h}^{0}}=h({{f}^{0}})$ of the optimal linear estimate of the functional ${{A}^{e}}\zeta $ is calculated by  formula (25) if the condition $h({{f}^{0}})\in {{H}_{{{D}_{f}}}}$ holds true.
\end{lem}
\begin{lem} \label{lem4.4}
 Spectral densities ${{f}^{0}}(\lambda )\in {{D}_{f}}$ and${{g}^{0}}(\lambda )\in {{D}_{g}}$ which satisfy the minimality condition (6), are least favorable in the class $D$ for the optimal estimation of the functional ${{A}^{f}}\zeta $, if the Fourier coefficients of matrix functions
$${{({{f}^{0}}(\lambda )+{{g}^{0}}(\lambda ))}^{-1}},\quad {{f}^{0}}(\lambda ){{({{f}^{0}}(\lambda )+{{g}^{0}}(\lambda ))}^{-1}},\quad {{f}^{0}}(\lambda ){{({{f}^{0}}(\lambda )+{{g}^{0}}(\lambda ))}^{-1}}{{g}^{0}}(\lambda )$$
determine matrices $U_{{}}^{0},\,V_{{}}^{0},\,W_{{}}^{0}$, which give a solution of the extremum problem
\begin{equation}\label{eq34}
\underset{(f,g)\in D}{\mathop{\max }}\,\left( \left\langle a,Wa \right\rangle +\left\langle {{U}^{-1}}Va,Va \right\rangle  \right)=\left\langle a,{{W}^{0}}a \right\rangle +\left\langle {{({{U}^{0}})}^{-1}}{{V}^{0}}a,{{V}^{0}}a \right\rangle .
\end{equation}
The minimax-robust spectral characteristic ${{h}^{0}}=h({{f}^{0}},{{g}^{0}})$ of the optimal linear estimate of the functional ${{A}^{f}}\zeta $ is calculated by formula (32) if the condition $h({{f}^{0}},{{g}^{0}})\in {{H}_{D}}$ holds true.
\end{lem}

\begin{zau} \label{rem4.1}
The least favorable spectral densities ${{f}^{0}}(\lambda )\in {{D}_{f}}$, ${{g}^{0}}(\lambda )\in {{D}_{g}}$ and the minimax-robust spectral characteristic $h({{f}^{0}},{{g}^{0}})\in {{H}_{D}}$ form a saddle point of the function $\Delta (h;f,g)$on the set ${{H}_{D}}\times D$.
The saddle point inequalities hold true if ${{h}^{0}}=h({{f}^{0}},{{g}^{0}})$, $h({{f}^{0}},{{g}^{0}})\in {{H}_{D}}$ and $({{f}^{0}},{{g}^{0}})$ is a solution to the conditional extremum problem
\begin{equation}\label{eq35}
\Delta (h({{f}^{0}},{{g}^{0}});f,g)\to \max,\quad (f,g)\in D,
\end{equation}
where the functional $\Delta (h({{f}^{0}},{{g}^{0}});f,g)$ is defined as
$$\Delta (h({{f}^{0}},{{g}^{0}});f,g)=\frac{1}{2\pi }\int_{-\pi }^{\pi }{\left[ {{A}^{\top }}({{e}^{i\lambda }}){{g}^{0}}(\lambda )+{{({{C}^{0}}({{e}^{i\lambda }}))}^{\top }} \right]}{{({{f}^{0}}(\lambda )+{{g}^{0}}(\lambda ))}^{-1}}f(\lambda )\times $$
$${{({{f}^{0}}(\lambda )+{{g}^{0}}(\lambda ))}^{-1}}\left[ {{g}^{0}}(\lambda )\overline{A({{e}^{i\lambda }})}+\overline{{{C}^{0}}({{e}^{i\lambda }})} \right]d\lambda +\frac{1}{2\pi }\int_{-\pi }^{\pi }{\left[ {{A}^{\top }}({{e}^{i\lambda }}){{f}^{0}}(\lambda )-{{({{C}^{0}}({{e}^{i\lambda }}))}^{\top }} \right]}\times $$
$${{({{f}^{0}}(\lambda )+{{g}^{0}}(\lambda ))}^{-1}}g(\lambda ){{({{f}^{0}}(\lambda )+{{g}^{0}}(\lambda ))}^{-1}}\left[ {{f}^{0}}(\lambda )\overline{A({{e}^{i\lambda }})}-\overline{{{C}^{0}}({{e}^{i\lambda }})} \right]d\lambda $$
for the extrapolation problem;
$$\Delta (h({{f}^{0}},{{g}^{0}});f,g)=\frac{1}{2\pi }\int_{-\pi }^{\pi }{\left[ {{A}_{N}}^{\top }({{e}^{i\lambda }}){{g}^{0}}(\lambda )+{{({{C}_{N}}^{0}({{e}^{i\lambda }}))}^{\top }} \right]}{{({{f}^{0}}(\lambda )+{{g}^{0}}(\lambda ))}^{-1}}f(\lambda )\times $$
$${{({{f}^{0}}(\lambda )+{{g}^{0}}(\lambda ))}^{-1}}\left[ {{g}^{0}}(\lambda )\overline{{{A}_{N}}({{e}^{i\lambda }})}+\overline{{{C}_{N}}^{0}({{e}^{i\lambda }})} \right]d\lambda +\frac{1}{2\pi }\int_{-\pi }^{\pi }{\left[ {{A}_{N}}^{\top }({{e}^{i\lambda }}){{f}^{0}}(\lambda )-{{({{C}_{N}}^{0}({{e}^{i\lambda }}))}^{\top }} \right]}\times $$
$${{({{f}^{0}}(\lambda )+{{g}^{0}}(\lambda ))}^{-1}}g(\lambda ){{({{f}^{0}}(\lambda )+{{g}^{0}}(\lambda ))}^{-1}}\left[ {{f}^{0}}(\lambda )\overline{{{A}_{N}}({{e}^{i\lambda }})}-\overline{{{C}_{N}}^{0}({{e}^{i\lambda }})} \right]d\lambda $$
for the interpolation problem;
$$\Delta (h({{f}^{0}},{{g}^{0}});f,g)=\frac{1}{2\pi }\int_{-\pi }^{\pi }{\left[ {{A}_{-}}^{\top }({{e}^{i\lambda }}){{g}^{0}}(\lambda )+{{({{D}^{0}}({{e}^{i\lambda }}))}^{\top }} \right]}{{({{f}^{0}}(\lambda )+{{g}^{0}}(\lambda ))}^{-1}}f(\lambda )\times $$
$${{({{f}^{0}}(\lambda )+{{g}^{0}}(\lambda ))}^{-1}}\left[ {{g}^{0}}(\lambda )\overline{{{A}_{-}}({{e}^{i\lambda }})}+\overline{{{D}^{0}}({{e}^{i\lambda }})} \right]d\lambda +\frac{1}{2\pi }\int_{-\pi }^{\pi }{\left[ {{A}_{-}}^{\top }({{e}^{i\lambda }}){{f}^{0}}(\lambda )-{{({{D}^{0}}({{e}^{i\lambda }}))}^{\top }} \right]}\times $$
$${{({{f}^{0}}(\lambda )+{{g}^{0}}(\lambda ))}^{-1}}g(\lambda ){{({{f}^{0}}(\lambda )+{{g}^{0}}(\lambda ))}^{-1}}\left[ {{f}^{0}}(\lambda )\overline{{{A}_{-}}({{e}^{i\lambda }})}-\overline{{{D}^{0}}({{e}^{i\lambda }})} \right]d\lambda $$
for the filtering problem.

In the case of interpolation (extrapolation) problem for the corresponding functional from observations of the process $\zeta (t)$ without noise the conditional extremum problem (35) can be rewritten as
\begin{equation}\label{eq36}
\Delta (h({{f}^{0}});f)\to \max,\quad f\in {{D}_{f}},
\end{equation}
where the functional $\Delta (h({{f}^{0}});f)$ is defined as
$$\Delta (h({{f}^{0}});f)=\frac{1}{2\pi }\int_{-\pi }^{\pi }{{{({{C}^{0}}({{e}^{i\lambda }}))}^{\top }}}{{({{f}^{0}}(\lambda ))}^{-1}}f(\lambda ){{({{f}^{0}}(\lambda ))}^{-1}}\overline{{{C}^{0}}({{e}^{i\lambda }})}\,d\lambda $$
for the extrapolation problem,
$$\Delta (h({{f}^{0}});f)=\frac{1}{2\pi }\int_{-\pi }^{\pi }{{{({{C}_{N}}^{0}({{e}^{i\lambda }}))}^{\top }}}{{({{f}^{0}}(\lambda ))}^{-1}}f(\lambda ){{({{f}^{0}}(\lambda ))}^{-1}}\overline{{{C}_{N}}^{0}({{e}^{i\lambda }})}\,d\lambda $$
for the interpolation problem.
\end{zau}

\section{{Minimax-robust spectral characteristics for given classes $D$}}

\subsection{Interpolation problem in the class $D_{M}^{-}$}

Consider the minimax estimation problem for the functional
$$A_{N}^{i}\zeta =\int_{0}^{(N+1)T}{a}(t)\zeta (t)dt$$
based on observations of the process $\zeta (t)$ at points of time $t\in \mathbb{R}\backslash [0,(N+1)T]$ under the condition that the spectral density $f(\lambda )$ of the generated vector stationary sequence $\{{{\zeta }_{j}},j\in \mathbb{Z}\}$ belongs to the class
$$D_{M}^{-}=\left\{ f(\lambda )\left| \frac{1}{2\pi }\int_{-\pi }^{\pi }{{{f}^{-1}}(\lambda )\,\cos (m\lambda )\,d\lambda =P(m),\,m=0,1,...,M} \right. \right\},$$
where the sequence $P(m)=\{{{p}_{kn}}(m)\}_{k,n=1}^{\infty },\,\,m=0,1,...,M,$ is such that $P(m)={{P}^{*}}(-m)$ and the matrix function $\sum\nolimits_{m=-M}^{M}{P(m){{e}^{im\lambda }}}$ is a nonnegative matrix with nonzero determinant. With the help of Lagrange multipliers method we can find that solution ${{f}^{0}}(\lambda )$ of the conditional extremum problem (36) satisfies the relation
$${{[{{({{f}^{0}}(\lambda ))}^{-1}}]}^{\top }}C_{N}^{0}({{e}^{i\lambda }}){{\left( C_{N}^{0}({{e}^{i\lambda }}) \right)}^{*}}{{[{{({{f}^{0}}(\lambda ))}^{-1}}]}^{\top }}=$$
$$={{[{{({{f}^{0}}(\lambda ))}^{-1}}]}^{\top }}\,\left( \sum\nolimits_{m=0}^{M}{{{{\vec{\alpha }}}_{m}}{{e}^{im\lambda }}} \right){{\left( \sum\nolimits_{m=0}^{M}{{{{\vec{\alpha }}}_{m}}{{e}^{im\lambda }}} \right)}^{*}}{{[{{({{f}^{0}}(\lambda ))}^{-1}}]}^{\top }},$$
where ${{\vec{\alpha }}_{m}},\,m=0,1,...,M,$ are Lagrange multipliers. The last equality holds if
$$\sum\nolimits_{j=0}^{N}{{{{\vec{c}}}_{j}}{{e}^{ij\lambda }}}=\sum\nolimits_{m=0}^{M}{{{{\vec{\alpha }}}_{m}}{{e}^{im\lambda }}}.$$

Consider the following cases:  $M\ge N$ and $M<N$.

Let $M\ge N$. In this case the Fourier coefficients of the matrix function ${{({{f}^{0}}(\lambda ))}^{-1}}$  determine the matrix $B_{N}^{0}.$ Thus, the extremum problem (36) is degenerate, and Lagrange multipliers ${{\vec{\alpha }}_{N+1}}=...={{\vec{\alpha }}_{M}}=\vec{0}$.
We take ${{\vec{\alpha }}_{N+1}}=...={{\vec{\alpha }}_{M}}=\vec{0}$, find ${{\vec{\alpha }}_{0}},...,{{\vec{\alpha }}_{N}}$ from the equation $B_{N}^{0}\alpha _{0}^{N}={{a}_{N}},$ where $\alpha _{0}^{N}={{({{\vec{\alpha }}_{0}},...,{{\vec{\alpha }}_{N}})}^{\top }}$ and come to conclusion that the least favorable density ${{f}^{0}}(\lambda )$ satisfies the relation
\begin{equation}\label{eq37}
{{f}^{0}}(\lambda )={{\left( \sum\nolimits_{u=-M}^{M}{P(m){{e}^{im\lambda }}} \right)}^{-1}}={{\left( \left( \sum\nolimits_{j=0}^{M}{{{A}_{j}}{{e}^{-ij\lambda }}} \right){{\left( \sum\nolimits_{j=0}^{M}{{{A}_{j}}{{e}^{-ij\lambda }}} \right)}^{*}} \right)}^{-1}}.
\end{equation}
So ${{f}^{0}}(\lambda )$ is spectral density of the vector autoregressive stochastic sequence of order $M$
\begin{equation}\label{eq38}
\sum\nolimits_{j=0}^{M}{{{A}_{j}}{{{\vec{\zeta }}}_{l-j}}}={{\vec{\varepsilon }}_{l}}.
\end{equation}
Let $M<N$. In this case the matrix $B_{N}^{{}}$ is determined by the Fourier coefficients of the matrix function ${{(f(\lambda ))}^{-1}}$. Matrices $P(m),\,m=0,...,M,$ are known and matrices $P(m),\,m=M+1,...,N,$ are unknown. The unknown Lagrange multipliers $\,{{\vec{\alpha }}_{m}},m=0,...,M,$ and matrices $P(m),\,m=M+1,...,N,$ can be found from the equation
	$${{B}_{N}}\alpha _{0}^{M}={{a}_{N}},\quad \alpha _{0}^{M}={{({{\vec{\alpha }}_{0}},...,{{\vec{\alpha }}_{M}},\vec{0},...,\vec{0})}^{\top }}.$$	
If the sequence of matrices $P(m),\,\,m=0,1,...,N,$ form a positive definite matrix function $\sum\nolimits_{m=-N}^{N}{P(m){{e}^{im\lambda }}}$ with nonzero determinant, the spectral density
\begin{equation}\label{eq39}
{{f}^{0}}(\lambda )={{\left( \sum\nolimits_{u=-N}^{N}{P(m){{e}^{im\lambda }}} \right)}^{-1}}={{\left( \left( \sum\nolimits_{j=0}^{N}{{{A}_{j}}{{e}^{-ij\lambda }}} \right){{\left(
\sum\nolimits_{j=0}^{N}{{{A}_{j}}{{e}^{-ij\lambda }}} \right)}^{*}} \right)}^{-1}}
\end{equation}
is least favorable and determines the vector autoregressive stochastic sequence of order $N$
\begin{equation}\label{eq40}
\sum\nolimits_{j=0}^{N}{{{A}_{j}}{{{\vec{\zeta }}}_{l-j}}}={{\vec{\varepsilon }}_{l}}.
\end{equation}
Thus our results can be summarized in the following statement. For more details see the article by Dubovetska and Moklyachuk (2012c).

\begin{thm} \label{thm5.1.1}
 Spectral density (37) of the vector autoregressive stochastic sequence (38) of order $M$is least favorable in the class $D_{M}^{-}$  for the optimal linear estimation of the functional $A_{N}^{i}\zeta $ if $M\ge N.$ In the case where $M<N$ and solutions $P(m),\,\,m=M+1,...,N,$ of the equation ${{B}_{N}}\alpha _{0}^{M}={{a}_{N}}$ with coefficients $P(m),\,\,m=0,1,...,M,$ form a positive definite matrix function $\sum\nolimits_{m=-N}^{N}{P(m){{e}^{im\lambda }}}$ with nonzero determinant, spectral density (39) of the vector autoregressive stochastic sequence (40) of order $N$ is  least favorable in $D_{M}^{-}.$
The minimax-robust spectral characteristic $h({{f}^{0}})$ is given by (14).
\end{thm}

\subsection{Minimax extrapolation in the class $D_{0}^{1}$}

Consider the minimax approach to the problem of estimation of the functional
$${{A}^{e}}\zeta =\int_{0}^{\infty }{a(t)\zeta (t)\,dt}$$
based on observations of the process $\zeta (t)$ at points of time $t<0$ under the condition that the spectral density $f(\lambda )$ of the generated vector stationary sequence $\{{{\zeta }_{j}},j\in \mathbb{Z}\}$ admits the canonical factorization (24) and belongs to the class
$$D_{0}^{1}=\left\{ f(\lambda )\left| \frac{1}{2\pi }\int_{-\pi }^{\pi }{f(\lambda )\,d\lambda =P} \right. \right\},$$
where $P=\{{{p}_{kn}}\}_{k,n=1}^{\infty }$ is a given nonnegative Hermitian matrix. According to remark 4.1, the least favorable spectral density in the class $D_{0}^{1}$  gives solution to the problem
\begin{equation}\label{eq41}
\Delta (h({{f}^{0}});f)=\int_{-\pi }^{\pi }{{{({{S}^{0}}({{e}^{i\lambda }}))}^{\top }}}{{Q}^{0}}(\lambda )f(\lambda ){{Q}^{0}}(\lambda )\overline{{{S}^{0}}({{e}^{i\lambda }})}\,\,d\lambda \to \max ,\,\,\,\,\,\,f\in D_{0}^{1}.
\end{equation}
With the help of Lagrange multipliers method we can find that solution ${{f}^{0}}(\lambda )$ of the conditional extremum problem (41) satisfies the relation
\begin{equation}\label{eq42}
{{[{{S}^{0}}({{e}^{i\lambda }})]}^{\top }}\overline{{{S}^{0}}({{e}^{i\lambda }})}={{[{{P}^{0}}(\lambda )]}^{\top }}\,\vec{\alpha }\,{{\vec{\alpha }}^{*}}\overline{{{P}^{0}}(\lambda )},
\end{equation}
where $\,\vec{\alpha }$ is the Lagrange multiplier, ${{S}^{0}}({{e}^{i\lambda }})=\sum\nolimits_{l=0}^{\infty }{{{(A{{d}^{0}})}_{l}}{{e}^{il\lambda }}},$ ${{(A{{d}^{0}})}_{l}}=\sum\nolimits_{j=l}^{\infty }{\vec{a}_{j}^{\top }{{d}^{0}}(j-l)},$ $l\ge 0,$${{P}^{0}}(\lambda )=\sum\nolimits_{u=0}^{\infty }{{{d}^{0}}(u){{e}^{-iu\lambda }}}.$  Relation (42) holds true if the sequence of matrices ${{d}^{0}}=\{{{d}^{0}}(u),\,\,u=0,1,...\}$ satisfies the system of equations
\begin{equation}\label{eq43}
\sum\nolimits_{p=o}^{\infty }{\sum\nolimits_{s=0}^{\infty }{\overline{{{{\vec{a}}}_{r+p}}}\vec{a}_{s+p}^{\top }d(s)}}=\overline{{\vec{\alpha }}}{{\vec{\alpha }}^{\top }}d(r),\,\,\,\,r=0,1,....
\end{equation}
Restrictions of the class $D_{0}^{1}$ lead to the following condition on the sequence ${{d}^{0}}(u),\,\,u=0,1,...$
\begin{equation}\label{eq44}
\sum\nolimits_{u=o}^{\infty }{d(u){{d}^{*}}(u)}=P.
\end{equation}
Then the following theorem holds true.

\begin{thm} \label{thm5.2.1} Spectral density
$${{f}^{0}}(\lambda )=\left( \sum\nolimits_{u=0}^{\infty }{{{d}^{0}}(u){{e}^{-iu\lambda }}} \right){{\left( \sum\nolimits_{u=0}^{\infty }{{{d}^{0}}(u){{e}^{-iu\lambda }}} \right)}^{*}}$$
of one-sided moving average sequence of the form (22) is least favorable in the class $D_{0}^{1}$ for the optimal linear estimation of the functional ${{A}^{e}}\zeta $. The sequence of matrices ${{d}^{0}}=\{{{d}^{0}}(u),\,\,u=0,1,...\}$ satisfies relations (43) and condition (44).
The minimax-robust spectral characteristic $h({{f}^{0}})$ is calculated by formula (25).
\end{thm}

\subsection{Filtering problem in the class $D_{0}^{2}\times {{D}_{\varepsilon }}$}

Consider the minimax estimation problem for the functional
$${{A}^{f}}\zeta =\int_{0}^{\infty }{a(t)\zeta (-t)\,dt}$$
based on observations of the process $\zeta (t)+\theta (t)$ at points of time $t\le 0$, under the condition that spectral densities $f(\lambda )$, $g(\lambda )$ of the generated vector stationary sequences $\{{{\zeta }_{j}},j\in \mathbb{Z}\}$, $\{{{\theta }_{j}},j\in \mathbb{Z}\}$  belong to classes
$$D_{0}^{2}=\left\{ f(\lambda )|\frac{1}{2\pi }\int_{-\pi }^{\pi }{Tr\,\,f(\lambda )\,d\lambda }={{P}_{\varsigma }} \right\},$$
$${{D}_{\varepsilon }}=\left\{ g(\lambda )|g(\lambda )=\varepsilon \cdot {{g}_{1}}(\lambda )+(1-\varepsilon )\cdot {{g}_{2}}(\lambda ),\,\,\frac{1}{2\pi }\int_{-\pi }^{\pi }{Tr\,\,g(\lambda )\,d\lambda }={{P}_{\theta }} \right\},$$
where ${{g}_{1}}(\lambda )\ge 0$ is an unknown function, and ${{g}_{2}}(\lambda )$ is a given function.

With the help Lagrange multipliers method we find that solution $({{f}^{0}},{{g}^{0}})$ of the conditional extremum problem (35) satisfies relations
\begin{equation}\label{eq45}
\left( {{g}^{0}}(\lambda )\overline{{{A}_{-}}({{e}^{i\lambda }})}+\overline{{{D}^{0}}({{e}^{i\lambda }})} \right)\left( A_{-}^{\top }({{e}^{i\lambda }}){{g}^{0}}(\lambda )+{{({{D}^{0}}({{e}^{i\lambda }}))}^{\top }} \right)={{\alpha }^{2}}{{\left( {{f}^{0}}(\lambda )+{{g}^{0}}(\lambda ) \right)}^{2}},
\end{equation}
\begin{equation}\label{eq44}
\left( {{f}^{0}}(\lambda )\overline{{{A}_{-}}({{e}^{i\lambda }})}-\overline{{{D}^{0}}({{e}^{i\lambda }})} \right)\left( A_{-}^{\top }({{e}^{i\lambda }}){{f}^{0}}(\lambda )-{{({{D}^{0}}({{e}^{i\lambda }}))}^{\top }} \right)=({{\beta }^{2}}+\varphi (\lambda )){{\left( {{f}^{0}}(\lambda )+{{g}^{0}}(\lambda ) \right)}^{2}},
\end{equation}
where $\alpha _{{}}^{2},\,\,\beta _{{}}^{2}$ are Lagrange multipliers, the function $\varphi (\lambda )\le 0$ and $\varphi (\lambda )=0$ if $Tr\,\,{{g}^{0}}(\lambda )\ge Tr\,\,(1-\varepsilon )\cdot {{g}_{2}}(\lambda )$.

\begin{thm} \label{thm5.3.1}
Let spectral densities ${{f}^{0}}(\lambda )\in D_{0}^{2},$$\,{{g}^{0}}(\lambda )\in {{D}_{\varepsilon }}$ satisfy the minimality condition (6).  Then spectral densities ${{f}^{0}}(\lambda ),\,\,{{g}^{0}}(\lambda )$ are least favorable in the class $D_{0}^{2}\times {{D}_{\varepsilon }}$ for the optimal linear estimation of the functional ${{A}^{f}}\zeta ,$ if they are determined by relations (45), (46), give solution to conditional extremum problem (34) and satisfy restrictions which determine the class $D_{0}^{2}\times {{D}_{\varepsilon }}$. The minimax-robust spectral characteristic $h({{f}^{0}},{{g}^{0}})$ of the estimate ${{\hat{A}}^{f}}\zeta $ is calculated by the formula (32). The value of the mean square error $\Delta ({{f}^{0}},{{g}^{0}})$ is calculated by formula (33).
\end{thm}

\section{{Conclusions}}

In this article we describe methods of solution of the problem of optimal linear estimation of functionals which depend on unknown values of periodically correlated (PC) processes. We present a transition procedure from PC processes $\{\zeta (t),\,t\in \mathbb{R}\}$ and $\{\theta (t),\,t\in \mathbb{R}\}$ to the corresponding generated vector stationary sequences $\{{{\zeta }_{j}},j\in \mathbb{Z}\}$ and $\{{{\theta }_{j}},j\in \mathbb{Z}\}$. Decomposition of stationary sequences $\{{{\zeta }_{j}},j\in \mathbb{Z}\}$ and $\{{{\theta }_{j}},j\in \mathbb{Z}\}$ with the help of a special basis in Hilbert space allows us to reduce the estimation problem for PC processes to the corresponding problem for stationary vector-valued sequences.

The Hilbert space projection method is exploited to finding optimal linear estimates of functionals $A_{N}^{i}\zeta ,$ ${{A}^{e}}\zeta ,$ $A_{N}^{e}\zeta ,$ ${{A}^{f}}\zeta $ based on observations of the PC process $\zeta (t)+\theta (t)$ with the PC noise process $\theta (t)$. Formulas for calculating mean square errors and spectral characteristics of the optimal linear estimates of the corresponding functionals are proposed in the case of spectral certainty. Formulas that determine the greatest value of mean square errors and the minimax estimates of functionals ${{A}^{e}}\zeta $, $A_{N}^{e}\zeta $  are presented. It is shown that the least favorable sequence for the optimal estimation of ${{A}^{e}}\zeta $ and $A_{N}^{e}\zeta $ is one-sided moving average stationary sequence generated by PC process from the class $\mathbf{Y}$. The minimax approach to the problem of estimation of linear functionals $A_{N}^{i}\zeta ,$ ${{A}^{e}}\zeta ,$ $A_{N}^{e}\zeta ,$ ${{A}^{f}}\zeta $ is analyzed in the case of spectral uncertainty for concrete classes $D$ of spectral density matrices. Least favorable spectral densities and minimax-robust spectral characteristics of the optimal estimates of functionals are determined.

\end{document}